\documentclass[final,reqno]{amsart}
\usepackage{pdfsync}
\usepackage{amsmath,amssymb}
\usepackage{enumerate}
\usepackage{xspace}
\usepackage{array}
\usepackage{tabularx}
\usepackage{booktabs}
\usepackage{subfigure}
\usepackage{caption}
\usepackage{cite}
\usepackage{caption}
\usepackage{james}
\usepackage{cleveref}
\usepackage{fullpage}
\AtBeginDocument{%
  \let\cref\Cref
}

\graphicspath{{figures/}}

\usepackage{tikz}
\usetikzlibrary{shapes,arrows,positioning,decorations.pathreplacing,
  intersections,matrix,patterns,fadings,calc,decorations.markings,arrows.meta}

\newlength\mylen

\let\oldnl\nl
\newcommand{\nonl}{\renewcommand{\nl}{\let\nl\oldnl}}

\hypersetup{colorlinks, linkcolor=blue, citecolor=blue, urlcolor=blue}

\newcommand{\revised}[1]{#1}

\usepackage{mathtools}
\numberwithin{equation}{section}

\author{James Jackaman}
\address{James Jackaman \thanks{ Department
    of Mathematical Sciences,
    NTNU, 7491 Trondheim, Norway {\tt{james.jackaman@ntnu.no}}.}}
\author{Scott MacLachlan}
\address{Scott MacLachlan \thanks{
    Department of Mathematics and Statistics, Memorial University of
    Newfoundland, St.\ John's, NL, A1C 5S7, Canada
    {\tt{smaclachlan@mun.ca}}.}}

\thanks{This work was partially supported by an NSERC Discovery Grant
  (SM) and the European Union’s Horizon 2020 research and innovation
  program under the Marie Sk{\l}odowska-Curie grant agreement No
  101108679 (JJ)}

\title{Space-time waveform relaxation multigrid for Navier-Stokes}
\date{\today}

\begin{document}

  \maketitle

  \begin{abstract}
  Space-time finite-element discretizations are well-developed in many
  areas of science and engineering, but much work remains within the
  development of specialized solvers for the resulting linear and
  nonlinear systems.  In this work, we consider the all-at-once
  solution of the discretized Navier-Stokes equations over a
  space-time domain using waveform relaxation multigrid methods.  In
  particular, we show how to extend the efficient spatial multigrid
  relaxation methods from~\cite{Rafiei_2024} to a waveform relaxation
  method, and \revise{provide a proof-of-concept of the algorithmic} efficiency of the resulting monolithic
  Newton-Krylov-multigrid solver.  Numerical results demonstrate the
  scalability of the solver for varying discretization order and
  physical parameters.
\end{abstract}


\section{Introduction} \label{sec:introduction}

High-fidelity computer simulation of fluid flow has been a common target in scientific computing since the 1960's~\cite{Chorin1968, FHHarlow_JEWelch_1965a}.  As both high-performance computing platforms and simulation technology have advanced, computational fluid dynamics has persisted as a driving application, seeking higher-fidelity models (e.g., by using pressure-robust finite-element methods~\cite{VJohn_etal_2017a}) and better use of highly parallel computing platforms (such as via GPU computing~\cite{fischer2022nekrs}).  On modern manycore HPC systems, however, it is easy to reach saturation with spatial parallelism alone, where wallclock time to solution plateaus as increased communication costs outweigh any speedups gained by increasing core counts.  Recent studies suggest this occurs for problem sizes below $\mathcal{O}(10,000)$ degrees of freedom per processor on many modern machines~\cite{10.1145/2938615.2938617}.

One way to achieve better performance in the strong scaling limit is the use of time-parallel algorithms (cf.~\cite{10.1007/978-3-319-23321-5_3}) to increase the number of degrees of freedom in concurrent computation.  Recent years have seen the development of several families of time-parallel algorithms for numerical integration of PDEs, including Parareal~\cite{parareal-2001}, multigrid reduction-in-time (MGRIT)~\cite{RDFalgout_etal_2014a}, ParaDIAG~\cite{gander2020paradiag}, waveform relaxation (both Schwarz and multigrid approaches)~\cite{ELelarasmee_etal_1982, CLubich_AOstermann_1987, MBjorhus_1995}, and others.  A common critique, however, is that these algorithms are often only demonstrated on simple discretizations of simple diffusive problems, and that their performance often suffers on more complex models (particularly hyperbolic PDEs) without substantial additional insights or costly stages in the algorithmic development~\cite{DeSterck_etal_2019b, FDanieli_SMacLachlan_2019a}.

The focus of this manuscript is the development of multigrid waveform relaxation techniques for the Navier-Stokes equations, modeling time-dependent incompressible fluid flow.  Waveform relaxation methods were first proposed in the 1980's for circuit simulation~\cite{ELelarasmee_etal_1982}.  Their adaptation as a relaxation scheme for space-time multigrid methods was studied extensively in the 1990's, following the first work by Lubich and Ostermann~\cite{CLubich_AOstermann_1987}.  Much of this work has focused on methods for the heat equation and similar scalar parabolic model problems~\cite{SVandewalle_GHorton_1995, JJanssen_SVandewalle_1994, SVandewalle_EVandeVelde_1994, JJanssen_SVandewalle_1996a, JJanssen_SVandewalle_1996b, JvanLent_SVandewalle_2002a}.  To our knowledge, the only previous application of multigrid waveform relaxation to the Navier-Stokes equations is the work of Oosterlee and Wesseling~\cite{CWOosterlee_PWesseling_1993}, but this work is limited to a low-order finite-volume discretization and BDF(2) time discretization.  In particular, since higher-order BDF methods are not A-stable~\cite{MR2657217}, extensions of this approach to higher-order discretizations are not immediately apparent.  Furthermore, the multigrid method in~\cite{CWOosterlee_PWesseling_1993} is based on a line relaxation strategy that is appropriate to the curvilinear coordinates used therein, but not needed for the simpler setting of Cartesian grids considered here.

The main contribution of this work is the extension of efficient
spatial relaxation schemes for finite-element discretizations of the
Poisson and time-steady Navier-Stokes equations to waveform relaxation
schemes for space-time finite-element discretizations.  We use
standard continuous Lagrange finite-element methods for the spatial
discretization of the heat equation, focusing on polynomial degrees
from 1 to 3, and making use of patch-based relaxation to obtain
convergence rates that are independent of polynomial order. For
Navier-Stokes, we use Taylor-Hood elements of both lowest and second-lowest order ($\cpoly{2}$--$\cpoly{1}$ and $\cpoly{3}$--$\cpoly{2}$), and demonstrate that the spatial relaxation proposed in~\cite{Rafiei_2024} also extends to an effective waveform relaxation scheme.  \revise{Similar algorithms have also been explored in the context of a single timestep of a higher-order temporal discretization, primarily for implicit Runge-Kutta discretizations~\cite{RAbuLabdeh_etal_2022a,doi:10.1137/23M1569344}.}  In both cases, we couple spatial semidiscretization with discontinuous-Galerkin timestepping, leading to space-time finite-element discretizations on prismatic grids.  We demonstrate that the resulting Newton-Krylov-multigrid solver is highly efficient across a range of discretization orders and grid sizes, for two standard model problems.  \revise{This work serves as motivation for an efficient parallel-in-time implementation of the underlying waveform relaxation scheme, since we demonstrate that it leads to scalable algorithms that could be extended using the cyclic reduction methodology~\cite{RDFalgout_etal_2015a}.  To further explore this, we consider a parallel performance model that shows expected speedups given sufficient computing resources.}

The remainder of this paper is structured as follows.  In~\cref{sec:fem}, we review the space-time finite-element formulation considered here.  The waveform relaxation multigrid algorithm is detailed in~\cref{sec:wrmg}, \revise{with a supporting parallel performance model developed in~\cref{sec:wrmg:performancemodel}.}  Numerical results for both the heat equation and the Navier-Stokes equations are presented in~\cref{sec:numerics}.  Discussion and concluding remarks follow in~\cref{sec:conclusions}.

\section{Space-time finite element discretizations} \label{sec:fem}

Throughout this work, we utilize space-time discretizations that, at
their core, rely on a tensor-product structure separating the treatment of spatial and temporal domains.
We make this choice as the PDEs under consideration possess
fundamentally different dynamics in space and time. As such, we
make use of tensor-product finite-element spaces and, for clarity, we
expose our temporal and spatial discretizations independently. Before
proceeding, we note that mixed space-time elements have proven
successful when developing adaptive strategies \cite{Moiola2017} and
are also well developed \cite{Varoglu1980, Hulbert1990, Steinbach2015}.

Here we consider upwind discontinuous Galerkin temporal elements, see
\cite{Estep2002}. This temporal discretization can be applied to any
initial value problem (after reformulating it as a first-order-in-time
system), so we consider the following: Let
$\vec{y} \in\bs{H^1([0,T],\shilbert)}^d$ where
$d\in\mathbb{Z}^+$, $t\in[0,T]$ and
$\vec{x}\in \Omega \subset \mathbb{R}^2$, and $\shilbert$ represents an
appropriate spatial Hilbert space. We first consider (evolutionary)
PDEs of the type
\begin{equation} \label{eqn:pde}
  \vec{y}_t + \mathcal{L} \vec{y} = 0
  \qquad
  \vec{y}(0,\vec{x}) = \vec{y}_0(\vec{x})
  ,
\end{equation}
subject to appropriate spatial boundary conditions, where
$\mathcal{L}$ is an
infinite dimensional (nonlinear) differential operator and
$\vec{y}_0\in\left(\shilbert\right)^d$ is some given initial
condition. Indeed, for any problem of this form, we can discretize
temporally in a unified manner.

\begin{definition}[Discontinuous Galerkin notation]
  
  Through partitioning our time interval into $N$ sub-intervals
  $0=t_0<t_1<\ldots<t_N=T$, which we write as $I_n:=(t_n,t_{n+1})$ for $0 \leq n < N$, we
  may define the discontinuous finite-element space of degree $q$
  as
  \begin{equation}
    \dpoly{q}([0,T],\shilbert)
    =
    \left\{ Y \in L^2([0,T],\shilbert) : Y \big|_{I_n} \in \mathbb{P}_q(I_n,\shilbert),
      \text{ for }n=0,\ldots,N-1 \right\}
    ,
  \end{equation}
  where $\mathbb{P}_{q}(I_n,\shilbert)$ represents the space of degree
  $q$ polynomials on $I_n$ temporally (and functions in $\shilbert$
  spatially). We observe that this function space does not prescribe
  values at the temporal nodes $t_n$, which leads us to introduce the
  shorthand notations
  \begin{equation}
    \vec{Y}^+_n := \lim_{t \nearrow t_n} \vec{Y}(t, \vec{x})
    \quad \textrm{and} \quad
    \vec{Y}^-_n := \lim_{t \searrow t_n} \vec{Y}(t, \vec{x})
    .
  \end{equation}
  For convenience, we further introduce the jump of a solution about a
  point $t_n$ as
  $
  \jump{\vec{Y}_n} := \vec{Y}_n^+ - \vec{Y}_n^-
  .
  $

\end{definition}

\begin{definition}[Temporal discretisation] \label{def:dgintime}
  
  Let $\vec{y}$ solve \eqref{eqn:pde}, then our finite element
  approximation is given by seeking
  $\vec{Y}(\cdot,\vec{x})\in \bs{\dpoly{q}([0,T],\shilbert)}^d$ such that
  \begin{equation} \label{eqn:dgintime}
      \sum_{n=0}^{N-1} \int_{I_n}
      \bc{ \vec{Y}_t + \mathcal{L}\vec{Y} }
      \cdot \vec{\psi} \di{t}
      + \sum_{n=0}^{N-1} \jump{\vec{Y}_n} \cdot \vec{\psi_n}^+ 
      =
      0
      \qquad
      \forall \vec{\psi} \in \bs{\dpoly{q}([0,T],\shilbert)}^d
      ,
    \end{equation}
    where we define $\vec{Y}_0^- := \vec{Y}_0 = \vec{y}_0(\vec{x})$ to be our initial
    condition in time.

\end{definition}

\begin{remark}[Efficient implementation in time]

  Typically, when exploiting a finite-element method in time of the
  form \eqref{eqn:dgintime}, one localizes the approximation to a
  single element and solves in a timestepping fashion. In this work,
  we primarily focus on solving \emph{globally} in time, motivated by
  the temporal parallelism of the waveform relaxation multigrid
  method. We compare against a more typical timestepping
  implementation of the method in~\cref{sec:numerics}.

\end{remark}
  
\begin{remark}[Conforming temporal discretisation]
  
  We note that for problems of the form \eqref{eqn:pde}, it is natural
  to seek $\vec{y} \in \bs{H^1([0,T],\shilbert)}^d$ such that
  \begin{equation}
    \int_0^T \bc{\vec{y}_t + \mathcal{L}\vec{y}} \psi \di{t}
    =
    0
    \qquad
    \forall \psi \in \bs{ L^2([0,T],\shilbert)}^d
    ,
  \end{equation}
  as this choice of (temporal) function spaces allows for the
  conditions of Lax-Milgram to be met, utilizing the fact that
  $\vec{y}_t\in \bs{L^2([0,T],\shilbert)}^d$. To accurately discretize the
  temporal dynamics, this formulation allows use of continuous
  Lagrange elements of degree $q$ for $\vec{Y}(\cdot,\vec{x})$ tested
  against discontinuous Lagrange elements of degree $q-1$. However,
  the conforming nature in time restricts the ability to adapt in
  space-time so, in practice, the non-conforming upwind discontinuous
  Galerkin approximation \eqref{eqn:dgintime} is often used. We note
  that this approximation adds numerical dissipation, the amount of
  which decreases as we increase the order of the
  method~\cite{Estep2002}.

\end{remark}

We will now proceed to discuss the problems under consideration and
their spatial discretization. Throughout, we exploit the spatial $L^2$
inner product
\begin{equation} \label{eqn:inner}
  \sinner{\vec{u}, \vec{v}}
  :=
  \int_{\Omega} \operatorname{tr}\bc{\vec{u}^T\vec{v} } \di{\vec{x}}
  .
\end{equation}
In the sequel, we shall abuse notation by allowing \eqref{eqn:inner}
to contract tensors to scalars (through the trace inner product).

We shall apply the algorithm developed in this work to two PDEs,
namely the two-dimensional heat equation and incompressible
Navier-Stokes equations.  While both the examples below are presented without external forcing, such terms are easily added to the resulting weak forms and discretizations.

\begin{definition}[Heat equation] \label{def:heat}

  Let $u=u(t,\vec{x})\in H^1\bc{[0,T],H^2(\Omega)}$ where
  $\Omega\subset\mathbb{R}^2$ solve the heat equation
  \begin{equation} \label{eqn:heat}
    u_t - \nabla^2 u = 0
    ,
  \end{equation}
  subject to initial data $u(0,\vec{x}) = u_0(\vec{x})$ and natural
  boundary conditions.  Further let $\mathcal{T}_h$ be a triangulation
  of the polygonal spatial domain and
  $\cpoly{k}(\mathcal{T}_h) \subset H^1(\Omega)$ be the space of
  continuous Lagrange elements of degree $k$. Then, recalling notation
  from our temporal discretization (Definition \ref{def:dgintime}),
  our full discretization is described by seeking
  $U\in\dpoly{q}\bc{[0,T],\cpoly{k}(\mathcal{T}_h)}$ such that
  \begin{equation} \label{eqn:heatfem}
    \sum_{n=0}^{N-1} \int_{I_n}
    \sinner{ U_t, \psi}
    \di{t}
    +
    \sum_{n=0}^{N-1} \jump{U_n} \psi_n^+
    +
    \int_0^T
    \sinner{\nabla U, \nabla \psi}
    \di{t}
    =
    0
    \qquad
    \forall\psi\in\dpoly{q}\bc{[0,T],\cpoly{k}(\mathcal{T}_h)}
    ,
  \end{equation}
  where we define $U_0^- := \rerevise{\mathcal{I}} u_0(\vec{x})$ to be
  \rerevise{the interpolant of some initial data $u_0(\vec{x})$ into
    $\cpoly{k}(\mathcal{T}_h)$}.

\end{definition}

\begin{definition}[Navier-Stokes equations]

      Let $\vec{u}=\vec{u}(t,\vec{x})\in\bs{H^1\bc{[0,T],H^2(\Omega)}}^2$
      and $\phi=\phi(t,\vec{x})\in L^2\bc{[0,T],H^1(\Omega)}$ where
      $\Omega\subset\mathbb{R}^2$ is some bounded spatial domain subject
      to appropriate boundary conditions. We consider the incompressible
      fluid model
      \begin{equation} \label{eqn:ns}
        \begin{split}
      \vec{u}_t + R \vec{u} \cdot \nabla \vec{u}
      + \nabla \phi - \nabla^2 \vec{u} & = 0 \\
      \nabla \cdot \vec{u} & = 0
                             ,
    \end{split}
  \end{equation}
  where $R$ is the Reynolds number\revise{, giving a non-dimensionalized ratio of intertial to viscous terms in the flow, typically defined as $R = uL/\nu$, where $u$ and $L$ are a characteristic flow velocity and lengthscale, respectively, and $\nu$ is the kinematic viscosity of the fluid}. We discretize spatially with
  Taylor-Hood elements, i.e., we seek
  $\vec{U}\in\bs{\dpoly{q}\bc{[0,T],\cpoly{k+1}(\Omega)}}^2$ and
  $\Phi \in \dpoly{q}\bc{[0,T],\cpoly{k}(\Omega)}$ such that
  \begin{equation} \label{eqn:nsfem}
    \begin{split}
      \sum_{n=0}^{N-1} \int_{I_n}
      \sinner{\vec{U}_t, \vec{\psi}} \di{t}
      +
      \sum_{n=0}^{N-1} \jump{\vec{U}_n} \cdot \vec{\psi}_n^+
      \qquad \qquad \qquad & \\
      +
      \int_0^T \sinner{R \vec{U} \cdot \nabla \vec{U}, \vec{\psi}}
      +
      \sinner{\Phi, \nabla \cdot \vec{\psi}}
      +
      \sinner{\nabla\vec{U},\nabla\vec{\psi}}
      \di{t}
      & =
        0
        \quad \forall \vec{\psi} \in \bs{\dpoly{q}\bc{[0,T],\cpoly{k+1}(\Omega)}}^2
      \\
      \int_0^T
      \sinner{\nabla \cdot \vec{U}, \chi}
      \di{t}
      & =
        0
        \quad \forall \chi \in \dpoly{q}\bc{[0,T],\cpoly{k}(\Omega)}
        .
    \end{split}
  \end{equation}
\end{definition}

We note that the Navier-Stokes equations do not precisely fall into the framework described above, since they form a system of differential algebraic equations (DAEs), but the discretization above arises from consistently applying the space-time discretization framework in this setting.
We shall consider two test cases for \eqref{eqn:nsfem} in this
work, both posed on domain $\Omega = [0,1] \times [0,1]$.

\begin{definition}[Chorin test problem~\cite{Chorin1968}] \label{def:chorin}

  Fixing zero Dirichlet boundary conditions on the normal component of $\vec{u}$, we obtain the
  Chorin test problem \cite{Chorin1968} which has exact solution given by velocity
  \begin{subequations} \label{eqn:chorin}
    \begin{equation}
      \begin{split}
        \vec{u}(t,\vec{x})
        =
        \begin{pmatrix}
          - \cos{\pi x_1} \sin{\pi x_2} \\
          \sin{\pi x_1} \cos{\pi x_2}
        \end{pmatrix} e^{-2\pi^2 t}
        ,
      \end{split}
    \end{equation}
    with corresponding pressure
    \begin{equation}
      \phi
      =
      R \frac{\pi}{4} \bc{
        \cos{2\pi x_1} + \cos{2\pi x_2}}
      e^{-4\pi^2 t}
      .
    \end{equation}
    \revised{The (weakly imposed) initial condition for $\vec{u}$ is chosen to match this expression at time $t=0$, while (strong) Dirichlet boundary conditions are imposed on $\vec{u}$ and are chosen to interpolate this expression restricted to $\partial\Omega$ at the temporal approximation points.  We note that an initial condition for $\phi$ is not needed in this discretization, since there is no time derivative of $\phi$ in the equations.}
  \end{subequations}

\end{definition}

\begin{remark}[Dependence on Reynolds number] \label{rem:reynolds}

  Typical solutions of the Navier-Stokes equations are expected to exhibit strong dependence on the Reynolds number, including steepening slopes of the velocity near domain boundaries where no-slip or no-flow boundary conditions are applied.  In contrast to this, the velocity solution to the Chorin test problem is independent of the Reynolds number, and the pressure solution only shows linear dependence.  Thus, while having an analytical solution for this test problem allows us to verify accuracy of our computed solutions, it does not provide a realistic test of solvers for Navier-Stokes in typical situations.  For this reason, we also consider the case of a lid-driven cavity, which demonstrates much more complex dynamics.
  
\end{remark}

\begin{definition}[Lid-driven cavity~\cite{Bozeman1973}] \label{def:lid}

  We fix zero velocity boundary conditions on both components of the velocity on the bottom, left and
  right boundaries, \revised{denoted collectively as $\Gamma_1 \subset \partial\Omega$,}
  \begin{equation}
    \left. \vec{u} \right|_{\Gamma_1}
    =
    \begin{pmatrix}
      0 \\
      0
    \end{pmatrix},
  \end{equation}
  and fix a velocity of
  \begin{equation}
    \left. \vec{u} \right|_{\Gamma_2}
    =
    \begin{pmatrix}
      1 \\
      0
    \end{pmatrix}
  \end{equation}
  on the top boundary\revised{, denoted by $\Gamma_2 \subset \partial\Omega$}. The test case is then initialized
  with zero velocity $\vec{u}$ and pressure $\phi$. Such an
  initialization lends itself, over time, to the well studied steady
  laminar solution, see, for example, \cite{Burggraf1966,
    Bozeman1973}.

\end{definition}

\section{Waveform relaxation multigrid} \label{sec:wrmg}

Classical multigrid methods make use of two complementary processes to damp errors in an approximation to the solution of a discretized system of differential equations, known as relaxation and coarse-grid correction.  In the setting of a simple elliptic (time-steady) PDE, relaxation is used to efficiently damp high-frequency errors in an approximation, leaving a smooth error, which can be effectively approximated by solution of a restricted problem on a coarser spatial grid.  While multigrid theory and practice were first developed for elliptic PDEs, questions on how to extend these approaches to time-dependent problems arose in the early 1980's~\cite{MR0806780}.  We note that applying spatial multigrid to a single timestep of a fully discretized time-dependent PDE is quite natural, particularly when implicit Euler or Crank-Nicolson type discretizations are used, and was considered even earlier in the development of multigrid~\cite{brandt1979multi}.  Instead, we are focusing on cases that solve the full space-time discretization over multiple timesteps concurrently, using relaxation and coarse-grid correction processes defined over all space-time grid points, and not just on a single time-line.

As noted above, the dynamic behavior of solutions to many time-dependent PDEs is dramatically different in the temporal direction than in the spatial direction.  As such, even the earliest works on multigrid for space-time problems typically considered coarsening only in the spatial direction, commonly known as semi-coarsening, an approach also considered for other problems with very anisotropic behavior.  While Hackbusch's original approach~\cite{MR0806780} simply used the space-time analogues of Jacobi or Gauss-Seidel point relaxation techniques, Lubich and Oostermann quickly proposed the use of so-called waveform relaxation~\cite{ELelarasmee_etal_1982} techniques as relaxation methods for space-time multigrid problems~\cite{CLubich_AOstermann_1987}.  This approach was extensively studied and developed in the 1990's, particularly by Vandewalle and collaborators~\cite{SVandewalle_GHorton_1995, JJanssen_SVandewalle_1994, SVandewalle_EVandeVelde_1994, JJanssen_SVandewalle_1996a, JJanssen_SVandewalle_1996b, JvanLent_SVandewalle_2002a}.

Here we consider a standard geometric multigrid interpretation of
waveform relaxation multigrid, \revised{as applied to a linear (or linearized) problem}. Given a hierarchy of spatial meshes,
$\{\mathcal{T}_H,\mathcal{T}_{H/2},\ldots,\mathcal{T}_h\}$, and fixed
temporal mesh $0 = t_0 < t_1 < \ldots < t_N = T$, we create a
hierarchy of space-time grids as
$\{t_n\}_{n=0}^N \otimes \mathcal{T}_H$,
$\{t_n\}_{n=0}^N\otimes \mathcal{T}_{H/2}$, $\ldots$,
$\{t_n\}_{n=0}^N \otimes \mathcal{T}_h$, where we use
Kronecker-product notation to denote the space-time mesh formed by a
tensor product of a temporal mesh with a spatial one.  We use standard
geometric multigrid V-cycles as preconditioners for FGMRES~\cite{YSaad_2003a}, where we
use FGMRES not for its flexibility (we use stationary V-cycles), but
because it explicitly stores the preconditioned Arnoldi vectors and,
thus, can reassemble the final solution without the additional
application of the preconditioner needed in classical
right-preconditioned GMRES.\revised{  Algorithm~\ref{alg:Vcycle} gives pseudocode for a standard V-cycle algorithm, highlighting the interleaving of relaxation and coarse-grid correction on every level of the hierarchy; for more details on multigrid methods in general (and alternative cycling strategies), see~\cite{WLBriggs_VEHenson_SFMcCormick_2000a, UTrottenberg_etal_2001a}.
\begin{algorithm2e}[!ht]
  \caption{Multigrid V-cycle}\label{alg:Vcycle}
  \For{grids $\mathcal{T}_h$, $\mathcal{T}_{2h}$, \ldots, $\mathcal{T}_{H/2}$}{%
    Pre-relaxation on current grid\\
    Compute current residual\\
    Restrict residual to next coarser grid, zero initial guess on next coarse grid
  }
  Direct solve on coarsest grid \\
  \For{grids $\mathcal{T}_{H/2}$, $\mathcal{T}_{H/4}$, \ldots, $\mathcal{T}_{h}$}{%
    Interpolate and add correction from next coarser grid\\
    Post-relaxation on current grid
    }
\end{algorithm2e}}

\revised{
To complete specification of the multigrid cycle, we must specify the grid-transfer operators to be used for restriction of residuals and interpolation of corrections, as well as the relaxation scheme to be used on each level for each problem.  Since the details are specific to the problem at hand, we discuss these sequentially.  One commonality, however, is that we use a Chebyshev-based acceleration of relaxation in both cases~\cite{https://doi.org/10.1002/nla.1979}.  In all experiments, we use V(2,2) cycles, where we apply two sweeps of the patch-based relaxation schemes described below as pre- and post-relaxation on every level.  As is the case for simple problems, this relaxation is most efficient when suitably weighted (i.e., when we choose appropriate step lengths for the correction computed by each relaxation sweep).  In this work, we use Chebyshev polynomials to determine the relaxation weights.}
To determine the interval that defines the
Chebyshev polynomial, we use \revised{10} steps of relaxation-preconditioned
GMRES to estimate the largest eigenvalue of the preconditioned system,
$\lambda$, and choose the interval to be $[0.25\lambda, 1.05\lambda]$,
where we choose the factor of $1/4$ because we are coarsening by a
factor of two in two spatial dimensions, and $1.05$ to provide a
safety factor, in case $\lambda$ is an underestimate of the true
largest eigenvalue.  These choices essentially give optimal relaxation
weights to reduce the errors associated with the upper $3/4$ of the
spectrum of the relaxation-preconditioned system (provided that
$\lambda$ is a good estimate of the largest eigenvalue of the
relaxation-preconditioned system).

We now describe the grid-transfer operators and waveform relaxation variants that we use for the two model problems considered here.  For the heat equation, we construct interpolation to space-time grid $\{t_n\}_{n=0}^N \otimes \mathcal{T}_h$ as $P = I \otimes P_h$, where $I$ is an identity operator over the temporal mesh and $P_h$ is the standard finite-element interpolation from grid $2h$ to grid $h$ for the $\bpoly{k}$ elements.  This expression for interpolation assumes that the degrees of freedom are ordered consistently with the tensor-product mesh expressions above but, in practice, this need not be the case.  For restriction from this grid, we use $R = P^T$.  \revised{That is, our interpolation and restriction operators are just the standard operators that would be used in geometric multigrid methods applied to the spatial discretization, ``copied'' over every approximation point in the temporal mesh.}  For relaxation for the heat equation, we use vertex-star relaxation~\cite{PCPATCH} in the spatial domain, extended (in the usual style of waveform relaxation) to couple all degrees of freedom in space-time that reside in the spatial patch.  Vertex-star patches are defined by taking each vertex in the mesh, identifying all edges and elements incident on the vertex, and then taking the set of DoFs associated with the vertex itself and those internal to the incident edges and elements.  A sample such patch for $\cpoly{2}$ on a regular triangular mesh is shown at left of Figure~\ref{fig:patches}.  \revised{These patches are then extended on the space-time mesh by coupling these spatial degrees of freedom across all approximation points in the temporal mesh.}  We note that such a spatial patch solve is known to lead to robust preconditioners for the Poisson equation (see, for example, \cite{Pavarino:1993}), so we expect to see robust performance for the heat equation using its waveform relaxation analogue.

For the Navier-Stokes equations, \revise{we consider Newton-Krylov-multigrid solvers that use inexact Newton's method with linear systems solved by waveform-relaxation multigrid preconditioned FGMRES with} two generalizations to the components used for the heat equation.  First, in usual (spatial) monolithic multigrid style, the spatial interpolation operator $P_h$ is now block diagonal,
\begin{equation}
P_h = \begin{bmatrix} P_{\vec{u}} & 0 \\ 0 & P_{\phi} \end{bmatrix},
\end{equation}
where $P_{\vec{u}}$ is the appropriate (spatial) finite-element
interpolation for the vector velocity space and $P_{\phi}$ is the
appropriate (spatial) finite-element interpolation for the scalar
pressure space.  Again, we use space-time interpolation operator $P = I \otimes P_h$ and restriction operator $R = P^T$, assuming compatible ordering of the discrete degrees of freedom\revise{, simply extending the spatial interpolation along each temporal approximation point}.  For relaxation, we use the waveform relaxation analogue of the Vanka+star relaxation recently proposed in~\cite{Rafiei_2024}.  Here, we define spatial patches around each vertex in the mesh.  Each patch contains all velocity degrees of freedom in the closure of the star of the vertex, which includes all velocity DoFs at the vertex itself, on the edges and elements incident on the vertex, and on edges and vertices directly adjacent to these edges and elements.  Along with these, we include all pressure DoFs on the vertex star.  \revise{Again, we extend these patches in the temporal direction to included these spatial DoFs at each temporal approximation point in one space-time patch.}  For a Taylor-Hood discretization with a $\cpoly{2}$ velocity space and $\cpoly{1}$ pressure space, this is the usual Vanka patch, but this patch has been observed to yield convergence that is robust in polynomial order for higher-order Taylor-Hood discretizations of time-independent problems.  A sample patch for the case of a $(\cpoly{3})^2$-$\cpoly{2}$ discretization is shown at right of Figure~\ref{fig:patches}.

\begin{figure}
          \centering
        \begin{tikzpicture}[scale=0.85]

        \draw[very thick] (0,0) -- (4,0) -- (4, 4) -- (0, 4) -- (0,0);
        \draw[very thick] (0,2) -- (4,2);
        \draw[very thick] (2,0) -- (2,4);
        \draw[very thick] (4,0) -- (0,4);
        \draw[very thick] (2,0) -- (0,2);
        \draw[very thick] (4,2) -- (2,4);
        \node[circle,draw,fill=gray!60,inner sep=0pt,minimum size=9pt] at (2,1) {};
        \node[circle,draw,fill=gray!60,inner sep=0pt,minimum size=9pt] at (3,1) {};
        \node[circle,draw,fill=gray!60,inner sep=0pt,minimum size=9pt] at (1,2) {};
        \node[circle,draw,fill=gray!60,inner sep=0pt,minimum size=9pt] at (2,2) {};
        \node[circle,draw,fill=gray!60,inner sep=0pt,minimum size=9pt] at (3,2) {};
        \node[circle,draw,fill=gray!60,inner sep=0pt,minimum size=9pt] at (1,3) {};
        \node[circle,draw,fill=gray!60,inner sep=0pt,minimum size=9pt] at (2,3) {};
        \node[circle,draw,fill=gray!10,inner sep=0pt,minimum size=9pt] at (0,0) {};
        \node[circle,draw,fill=gray!10,inner sep=0pt,minimum size=9pt] at (1,0) {};
        \node[circle,draw,fill=gray!10,inner sep=0pt,minimum size=9pt] at (2,0) {};
        \node[circle,draw,fill=gray!10,inner sep=0pt,minimum size=9pt] at (3,0) {};
        \node[circle,draw,fill=gray!10,inner sep=0pt,minimum size=9pt] at (4,0) {};
        \node[circle,draw,fill=gray!10,inner sep=0pt,minimum size=9pt] at (0,1) {};
        \node[circle,draw,fill=gray!10,inner sep=0pt,minimum size=9pt] at (1,1) {};
        \node[circle,draw,fill=gray!10,inner sep=0pt,minimum size=9pt] at (4,1) {};
        \node[circle,draw,fill=gray!10,inner sep=0pt,minimum size=9pt] at (0,2) {};
        \node[circle,draw,fill=gray!10,inner sep=0pt,minimum size=9pt] at (4,2) {};
        \node[circle,draw,fill=gray!10,inner sep=0pt,minimum size=9pt] at (0,3) {};
        \node[circle,draw,fill=gray!10,inner sep=0pt,minimum size=9pt] at (3,3) {};
        \node[circle,draw,fill=gray!10,inner sep=0pt,minimum size=9pt] at (4,3) {};
        \node[circle,draw,fill=gray!10,inner sep=0pt,minimum size=9pt] at (0,4) {};
        \node[circle,draw,fill=gray!10,inner sep=0pt,minimum size=9pt] at (1,4) {};
        \node[circle,draw,fill=gray!10,inner sep=0pt,minimum size=9pt] at (2,4) {};
        \node[circle,draw,fill=gray!10,inner sep=0pt,minimum size=9pt] at (3,4) {};
        \node[circle,draw,fill=gray!10,inner sep=0pt,minimum size=9pt] at (4,4) {};
        \node at (1.6,1.6) {};

        \tikzset{->-/.style={decoration={markings,mark=at position .6 with {\arrow[scale=1.7]{stealth}}},postaction={decorate}}}

        \foreach \x in {9,...,12}{
                \foreach \y in {0,...,3}{
                        \node[circle,draw,fill=gray!60,inner sep=0pt,minimum size=9pt] at (2*\x/3,2*\y/3+2) {};
                }
        }
        \foreach \x in {12,...,15}{
                \foreach \y in {0,...,3}{
                        \node[circle,draw,fill=gray!60,inner sep=0pt,minimum size=9pt] at (2*\x/3,2*\y/3) {};
                }
        }
        \node[circle,draw,fill=gray!60,inner sep=0pt,minimum size=9pt] at (22/3,2/3) {};
        \node[circle,draw,fill=gray!60,inner sep=0pt,minimum size=9pt] at (22/3,4/3) {};
        \node[circle,draw,fill=gray!60,inner sep=0pt,minimum size=9pt] at (20/3,4/3) {};
        \node[circle,draw,fill=gray!60,inner sep=0pt,minimum size=9pt] at (26/3,8/3) {};
        \node[circle,draw,fill=gray!60,inner sep=0pt,minimum size=9pt] at (26/3,10/3) {};
        \node[circle,draw,fill=gray!60,inner sep=0pt,minimum size=9pt] at (28/3,8/3) {};
        \draw[very thick] (6,0) -- (8,0);
        \draw[very thick] (8,0) -- (10,0);
        \draw[very thick] (6,2) -- (8,2);
        \draw[very thick] (8,2) -- (10,2);
        \draw[very thick] (6,4) -- (8,4);
        \draw[very thick] (8,4) -- (10,4);
        \draw[very thick] (6,0) -- (6,2);
        \draw[very thick] (6,2) -- (6,4);
        \draw[very thick] (8,0) -- (8,2);
        \draw[very thick] (8,2) -- (8,4);
        \draw[very thick] (10,0) -- (10,2);
        \draw[very thick] (10,2) -- (10,4);
        \draw[very thick] (8,0) -- (6,2);
        \draw[very thick] (10,0) -- (8,2);
        \draw[very thick] (8,2) -- (6,4);
        \draw[very thick] (10,2) -- (8,4);
        \node[rectangle,fill=black,inner sep=0pt,minimum size=4pt] at (8,2) {};
        \node[rectangle,fill=black,inner sep=0pt,minimum size=4pt] at (7,2) {};
        \node[rectangle,fill=black,inner sep=0pt,minimum size=4pt] at (9,2) {};
        \node[rectangle,fill=black,inner sep=0pt,minimum size=4pt] at (8,1) {};
        \node[rectangle,fill=black,inner sep=0pt,minimum size=4pt] at (9,1) {};
        \node[rectangle,fill=black,inner sep=0pt,minimum size=4pt] at (7,3) {};
        \node[rectangle,fill=black,inner sep=0pt,minimum size=4pt] at (8,3) {};
        \node[circle,draw,fill=gray!10,inner sep=0pt,minimum size=9pt] at (18/3,0) {};
        \node[circle,draw,fill=gray!10,inner sep=0pt,minimum size=9pt] at (18/3,2/3) {};
        \node[circle,draw,fill=gray!10,inner sep=0pt,minimum size=9pt] at (18/3,4/3) {};
        \node[circle,draw,fill=gray!10,inner sep=0pt,minimum size=9pt] at (20/3,0) {};
        \node[circle,draw,fill=gray!10,inner sep=0pt,minimum size=9pt] at (20/3,2/3) {};
        \node[circle,draw,fill=gray!10,inner sep=0pt,minimum size=9pt] at (22/3,0) {};
        \node[circle,draw,fill=gray!10,inner sep=0pt,minimum size=9pt] at (30/3,12/3) {};
        \node[circle,draw,fill=gray!10,inner sep=0pt,minimum size=9pt] at (30/3,10/3) {};
        \node[circle,draw,fill=gray!10,inner sep=0pt,minimum size=9pt] at (30/3,8/3) {};
        \node[circle,draw,fill=gray!10,inner sep=0pt,minimum size=9pt] at (28/3,12/3) {};
        \node[circle,draw,fill=gray!10,inner sep=0pt,minimum size=9pt] at (28/3,10/3) {};
        \node[circle,draw,fill=gray!10,inner sep=0pt,minimum size=9pt] at (26/3,12/3) {};
        \foreach \x in {6,...,10}{
          \node[rectangle,draw,fill=gray!30,inner sep=0pt,minimum size=4pt] at (\x,0) {};
          \node[rectangle,draw,fill=gray!30,inner sep=0pt,minimum size=4pt] at (\x,4) {};
          }
        \node[rectangle,draw,fill=gray!30,inner sep=0pt,minimum size=4pt] at (6,1) {};
        \node[rectangle,draw,fill=gray!30,inner sep=0pt,minimum size=4pt] at (7,1) {};
        \node[rectangle,draw,fill=gray!30,inner sep=0pt,minimum size=4pt] at (10,1) {};
        \node[rectangle,draw,fill=gray!30,inner sep=0pt,minimum size=4pt] at (6,2) {};
        \node[rectangle,draw,fill=gray!30,inner sep=0pt,minimum size=4pt] at (10,2) {};
        \node[rectangle,draw,fill=gray!30,inner sep=0pt,minimum size=4pt] at (6,3) {};
        \node[rectangle,draw,fill=gray!30,inner sep=0pt,minimum size=4pt] at (9,3) {};
        \node[rectangle,draw,fill=gray!30,inner sep=0pt,minimum size=4pt] at (10,3) {};
        \node at (7.6,1.6) {} ;

        \node[circle,draw,fill=gray!60,inner sep=0pt,minimum size=9pt] at (-0.5,-1) {};
        \node[right=3mm] at (-0.5,-1) {included $\cpoly{2}$ heat DoFs};
        \node[circle,draw,fill=gray!60,inner sep=0pt,minimum size=9pt] at (5.5,-1) {};
        \node[right=3mm] at (5.5,-1) {included $\cpoly{3}$ velocity DoFs};
        \node[circle,draw,fill=gray!10,inner sep=0pt,minimum size=9pt] at (5.5,-1.7) {};
        \node[right=3mm] at (5.5,-1.7) {excluded $\cpoly{3}$ velocity DoFs};
        \node[circle,draw,fill=gray!10,inner sep=0pt,minimum size=9pt] at (-0.5,-1.7) {};
        \node[right=3mm] at (-0.5,-1.7) {excluded $\cpoly{2}$ heat DoFs};
        \node[rectangle,fill=black,inner sep=0pt,minimum size=4pt] at (5.5,-2.4) {};
        \node[right=3mm] at (5.5,-2.4) {included $\cpoly{2}$ pressure DoFs};
        \node[rectangle,draw,fill=gray!30,inner sep=0pt,minimum size=4pt] at (5.5,-3.1) {};
        \node[right=3mm] at (5.5,-3.1) {excluded $\cpoly{2}$ pressure DoFs};

        \end{tikzpicture}
\caption{At left, a sample spatial vertex-star patch used for relaxation for the $\bpoly{2}$ elements for the heat equation.  At right, a sample Vanka+star patch for the $(\cpoly{3})^2$-$\cpoly{2}$ discretization of the Navier-Stokes equations.}\label{fig:patches}
\end{figure}

\revise{
We note three important differences between our algorithm and the only previous application of waveform relaxation to the Navier-Stokes equations in~\cite{CWOosterlee_PWesseling_1993}.  First, we consider quite a different set of discretizations, with~\cite{CWOosterlee_PWesseling_1993} considering a low-order finite-volume in space coupled with $\theta$-scheme or BDF(2) time discretizations, which leads to much simpler coupling between timesteps than occurs with fully implicit Runge-Kutta or Galerkin-in-Time discretizations.  The naturally leads to differences in the spatial interpolation and restriction operators.  Secondly, we consider a standard Newton-Krylov-multigrid approach to handling the nonlinear Navier-Stokes equations, while~\cite{CWOosterlee_PWesseling_1993} uses a nonlinear multigrid scheme.  Even for simple problems, our experience is that Newton-Krylov-multigrid tends to be more efficient than nonlinear multigrid (due to the high cost of nonlinear relaxation schemes), although comparing these two approaches may be of interest for future work.  Finally, \cite{CWOosterlee_PWesseling_1993} uses a line-relaxation version of patch-based relaxation on their finite-volume grids.  While this appears appealing due to its robustness, we did not find the need to use such coupling here, and it certainly increases the cost of the relaxation sweeps to do so.}

We implement these relaxation schemes in Firedrake~\cite{FiredrakeUserManual}, taking advantage of the ``extruded mesh'' functionality~\cite{gmd-9-3803-2016}.  To do so, we define a spatial mesh hierarchy by defining the coarsest grid, $\mathcal{T}_H$, and uniformly refining it a given number of times.  Then, we create a hierarchy of space-time grids as an extruded mesh hierarchy with a fixed (uniform) temporal grid at each level.  This allows us to specify space-time function spaces for the discretization by defining spatial and temporal finite elements, on triangles and intervals, respectively, then defining tensor-product elements, and defining function spaces to use these tensor-product elements over the space-time grids.  The extruded mesh construction lets us easily construct patches for waveform relaxation, as we can do this using the PETSc DMPlex abstraction to manage the construction~\cite{doi:10.1137/15M1026092}.  As extruded meshes are implemented in Firedrake, all space-time DoFs are assigned to spatial topological entities based on their spatial positioning on the grid, so that we can define patches for space-time extruded hierarchies in the same way as we do for spatial hierarchies, just by looping over the topological entities in the DMPlex structure.

\rerevise{Results in \cref{sec:numerics} compare the WRMG solver to its timestepping analogue, for which we use monolithic multigrid on the coupled equations for a single element in time.  We implement this using the Irksome framework~\cite{farrell2021irksome, RCKirby_SPMacLachlan_2024a}, using its extension to Galerkin-in-time methods from~\cite{BDAndrews_2026a}.  We note that the theory and practice of linear and nonlinear solvers for these discretizations is quite unstudied until now, but follow the practice of~\cite{doi:10.1137/23M1569344, RAbuLabdeh_etal_2022a} from the implicit Runge-Kutta case in the design of these solvers.  Mathematically, nothing changes from the WRMG solvers specified above, we simply consider couplings across only one time element in each timestep.}

\subsection{Parallel performance modelling}
\label{sec:wrmg:performancemodel}

\revise{
  Due to current software limitations, we are unable to perform parallel-in-time numerical experiments. This leads us to consider parallel performance models~\cite{Rauber_Runger_2000, 1592815, 10.1145/1995896.1995924} to predict the possible speedups that could be obtained on large-scale parallel machines with a performant implementation of WRMG.  We note that, as is almost always the case for parallel-in-time methods, temporal parallelism is necessary to achieve speedup, as the WRMG algorithm requires ``extra'' work compared to classical timestepping, but that work can be performed effectively in parallel.

  We develop two parallel performance models in this work, one for the cost of spatial multigrid as a solver in a classical space-parallel, but sequential-in-time, (timestepping) approach and one for the WRMG solver detailed above.  At the core of both of these solvers are patch solves, over a single temporal element in the timestepping case~\cite{RAbuLabdeh_etal_2022a,doi:10.1137/23M1569344} and over a complete timeline in the waveform relaxation case.  A natural question is how much these solves cost, since the patch matrices over each temporal element are somewhat sparse, but not substantially so.  To answer this question, we first note that the waveform relaxation patch solves have a block lower bidiagonal structure determined by the temporal elements, with the degrees of freedom on each temporal element only depending on those on the preceeding (in time) element.  Thus, the space-time patch solves for the waveform relaxation involve matrices of the form
  \begin{equation}\label{eq:waveform_patch_solve}
    \begin{bmatrix} A_{1,1} & \\
  A_{2,1} & A_{2,2} & \\
  & A_{3,2} & A_{3,3} & \\
  & & \ddots & \ddots & \\
  & & & A_{N,N-1} & A_{N,N} \end{bmatrix}
\begin{bmatrix} \vec{y}_1 \\ \vec{y}_2 \\ \vec{y}_3 \\ \vdots \\ \vec{y}_N \end{bmatrix}
= \begin{bmatrix} \vec{b}_1 \\ \vec{b}_2 \\ \vec{b}_3 \\ \vdots \\ \vec{b}_N \end{bmatrix},
\end{equation}
where the diagonal blocks are of dimension equal to the spatial patch size times the number of degrees of freedom on a single element in time ($q$+1 for discontinuous elements of degree $q$) and the subdiagonal blocks have rank equal to the number of spatial DoFs in the patch, coming from the jump term in~\eqref{eqn:dgintime}.  Thus, the cost of solving a single waveform relaxation sweep is roughly equal to $N$ times the cost of solving a single patch system over a single element in time.

To model the cost of a full solve, we first generated ``flame graphs'' showing the typical costs within a full run of our simulation code.  From this, we noted two main contributions to the program runtime.
One is the cost associated with assembling the finite-element interpolation operator for the high-order space-time systems.  For this, we model the cost as a scalar rate, $c_1$, times the total number of DoFs over the spatial grid and a single temporal element on each core (spatial DoFs per core times the number of temporal DoFs per element) raised to an unknown power, $p_1$, since this interpolation operator can be assembled for one element and used on all time-steps:
\[
T_{\text{interpolation assembly}} = c_1 \left((q+1)(\text{spatial DoFs per core})\right)^{p_1}.
\]
The other main cost is the (many) calls to the relaxation scheme.  To account for this cost, we propose a model for the cost of a single patch solve for the space-time system as equalling a scaling parameter, $c_2$, (effectively a FLOPS rate) times the number of timesteps, $N$, times the number of temporal DoFs on an element, $q+1$, raised to an unknown power, $p_2$, times the number of nonzeros in a spatial patch matrix to another unknown power, $p_3$:
\[
T_{\text{relaxation per patch}} = c_2 N (q+1)^{p_2}(\text{spatial nonzeros per patch})^{p_3}.
\]
We note that assuming the cost depends on the number of nonzero entries in a diagonal\rerevise{-block} matrix, $A_{i,i}$, seems reasonable, but that this is reflected quite differently between the temporal degrees of freedom (that ``stack'' small dense blocks into the matrix) and the spatial ones, that have more complex connections.  
For a single V-cycle with $M_x \times M_x$ spatial vertices on the finest grid, with regular coarsening by a factor of 2 in each direction on each level, we expect to perform 4 sweeps of relaxation on each level (since we use V(2,2) cycles), giving
\[
4\left( M_x^2 + (M_x/2)^2 + (M_x/4)^2 + ... \right) \approx \frac{16}{3} M_x^2
\]
patches for relaxation split evenly across $p$ cores.  This gives us a total cost model of
\begin{align}
T_{\text{solve}} =  & c_1 \left((q+1)(\text{spatial DoFs per core})\right)^{p_1} \label{eq:solve_cost_model} \\ & + c_2 \frac{16 V N M_x^2}{3p}(q+1)^{p_2}(\text{spatial nonzeros per patch})^{p_3}, \notag
\end{align}
where $V$ is the number of multigrid V-cycles needed for the solve.

Counting the spatial DoFs and nonzeros is tedious but straightforward.  For the heat equation discretized with order $k$ continuous elements on an $M_x \times M_x$ grid, we can simply count that we have 1 DoF per vertex, $k-1$ DoFs per edge, and $(k-1)(k-2)/2$ DoFs per element~\cite{Kirby2012}, giving
\begin{equation}
  \label{eq:PkDoFcount}
M_x^2 + (3(M_x-1)^2 + 2(M_x-1))(k-1) + 2(M_x-1)^2(k-1)(k-2)/2
\end{equation}
spatial DoFs, where the edge terms account for $M_x(M_x-1)$ edges aligned with each of the coordinate axes and $(M_x-1)^2$ diagonal edges in the mesh, and we count $2(M_x-1)^2$ triangular elements.  Similarly, each vertex-star patch (as shown in Figure~\ref{fig:patches}) has one vertex, 6 edges, and 6 elements, so each patch has
\begin{equation}
  \label{eq:PkPatchDoFs}
1 + 6(k-1) + 6(k-1)(k-2)/2 \rerevise{= 3k^2 - 3k + 1}
\end{equation}
spatial DoFs.  To count the nonzeros in the patch matrix, we count that the central vertex has connections to all DoFs in the patch, while each edge DoF has connections to the central vertex DoF, all of the DoFs on 3 edges (its own and 2 ``adjacent'' edges), and all of the DoFs on the two adjacent elements.  Similarly, each (internal) elemental DoF has connections to the central vertex DoF, all of the DoFs on the two adjacent edges in the patch, and all of the DoFs on its own element.  This leads to a count of
\begin{align*}
  1 + 6(k-1) + 6(k-1)(k-2)/2 & + 6(k-1)\left(3k-2 + 2(k-1)(k-2)/2\right)\\ & + 6((k-1)(k-2)/2)((k+1)(k+2)/2) \\
  & \rerevise{= \frac{3}{2}k^4 + 6k^3 - \frac{21}{2}k^2 - 3k + 7}
\end{align*}
\rerevise{nonzeros in the patch matrix}.

Similar (but more complex) counting is possible for the Navier-Stokes equations, where the spatial discretization uses $\cpoly{k+1}$ elements for each of the two components of the velocity, and $\cpoly{k}$ elements for the pressure.  Following from~\eqref{eq:PkDoFcount}, this gives
\[
\rerevise{3(M_x-1)^2k^2 +\left(4(M_x-1)^2+6(M_x-1)\right)k +2(M_x-1)^2+4(M_x-1)+3}
\]
DoFs over the entire spatial mesh.  Counting DoFs in the patches, we have the same number of pressure DoFs as in~\eqref{eq:PkPatchDoFs}, but the velocity DoFs cover a larger region, containing an additional 6 vertices and edges, in addition to having 2 components and being of order $k+1$.  Accounting for this, we have
\[
\rerevise{9k^2 + 15k + 15}
\]
DoFs in a single Vanka-star patch for this discretization.  These additional DoFs must be carefully accounted for in the nonzero entries in the patch system, along with separately counting the velocity-velocity nonzeros (which appear only once in the patch system) and the velocity-pressure nonzeros (that appear twice, once in the (1,2) block of the patch matrix and once in the (2,1) \rerevise{block}).  There are no pressure-pressure nonzeros in the patch.  \rerevise{Details of this calculation are given in~\cref{ssec:patch_nonzeros}, resulting in a count of
\[
12k^4 + 96k^3 + 240k^2 + 252k + 128
\]
nonzeros in a typical patch matrix.
We note that while the number of DoFs in a patch scales like $k^2$ and there are $\mathcal{O}(k^4)$ nonzero elements in each patch matrix, there are also many zero entries.  Thus, in practice, we use a sparse direct solver to factorize the patch matrices to take advantage of what sparsity does appear. }

To estimate the constants in~\eqref{eq:solve_cost_model}, we take the measured times-to-solution for the heat equation in Tables~\ref{tab:heat:time3} and~\ref{tab:heat:time4} below with $N = 20$, and use the least squares solver from \texttt{scipy.optimize} to minimize the least squares error between the model cost and the measured time.  This finds a relatively good fit, with less than 10\% error for all but one data point with $M_x = 161$, but slightly larger errors for the smaller grid with $M_x = 81$, where the times are smaller.  The constants we find are $c_1 = 5.9\times 10^{-4}$ and $c_2 = 3.5\times 10^{-7}$, and $p_1 = 1.29$, $p_2 = 2.45$, and $p_3 = 1$.  We note that $c_2$ has a more direct interpretation as a FLOPS rate (since we carefully count operations for the patch solves), while $c_1$ accounts for many operations in the assembly of the space-time interpolation operator. \rerevise{We plot this data fit in Figure \ref{fig:pmvalidation} for the data in Tables~\ref{tab:heat:time3} and~\ref{tab:heat:time4}.} We use these parameters in our model to predict solver costs for both the classical time-stepping algorithm and the \rerevise{time-parallelized} waveform relaxation approach as applied to the Navier-Stokes equations, where memory limitations prevented us from experimenting over a broader range of discretization orders.

\begin{figure}[h!]
  \centering
  \includegraphics[width=0.67\textwidth]{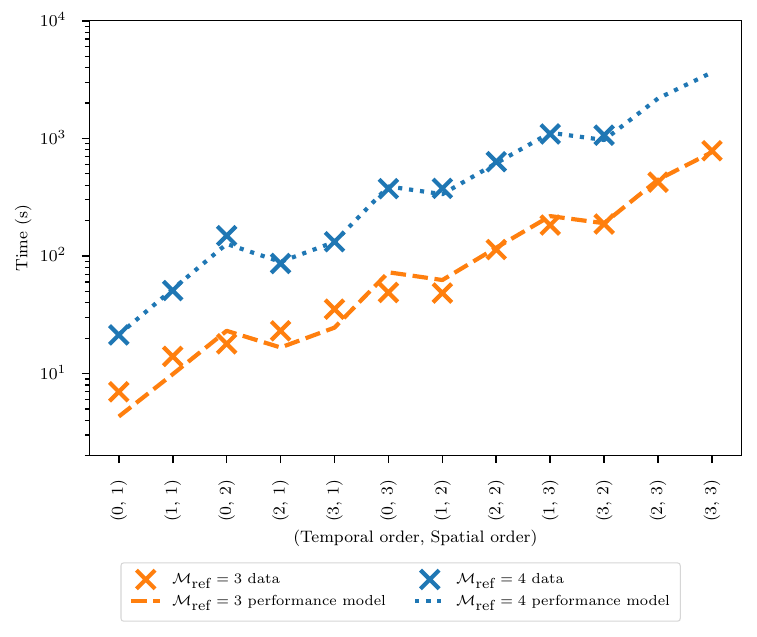}

  \caption{\rerevise{ Comparison of the measured time-to-solution for the heat equation with the model in~\eqref{eq:solve_cost_model} for multiple levels of spatial refinement
      ($\Mref$), and varying spatial and temporal degrees.  Measured time-to-solution for the heat equation (data from Tables~\ref{tab:heat:time3}
      and~\ref{tab:heat:time4}) is marked as crosses. The dotted or
      dashed lines outline the fit for the performance model in~\eqref{eq:solve_cost_model}  }  \label{fig:pmvalidation}}
  
\end{figure}

To predict performance for WRMG, we use a standard model of parallel communication time, modelled by a communication cost for a message of $n$ values (double-precision real numbers) of
\[
    T_{\text{comm}} = \alpha + n\beta,
  \]
  where $\alpha$ is the network latency (in seconds) and $\beta$ is the inverse bandwidth (in seconds per value), and a computational cost for $n$ floating-point operations of
  \[
    T_{\text{comp}} = n\gamma.
  \]
  We adapt the computational cost model in~\eqref{eq:solve_cost_model} for $T_{\text{comp}}$ and discuss the details of the communication model \rerevise{in~\cref{ssec:halos}.  We} assume we have an $M_x \times M_x$ finest grid that we parallelize evenly across $p_x \times p_x$ processors.  For the WRMG case, we assume that we evenly parallelize $N$ timesteps across $p_t$ processors.  We also assume that we will coarsen down to $\mathcal{O}(1)$ degrees of freedom on the spatial mesh, so we have $\logtwo{M_x}$ levels in the multigrid hierarchy.

\rerevise{For spatial multigrid (timestepping), the computational cost model retains} the first term in~\eqref{eq:solve_cost_model} (since this needs to be done once in the spatial multigrid case as well), \rerevise{while} the second term only \rerevise{replaces} $p$ in the denominator by $p_x^2$, the number of processors assumed. Note that this ``exchanges'' simultaneous work on $N$ timesteps within WRMG in that model with $N$ \rerevise{sequential} timesteps, each solved using V-cycles.  

\rerevise{
For WRMG, we need to consider parallel solution of the lower bidiagonal linear system from~\eqref{eq:waveform_patch_solve}.}  Using cyclic reduction (see~\cite{RDFalgout_etal_2015a, REISNER2020102705}), we can solve the space-time patch system by first computing a ``decoupled'' system for the first time step of each patch on each processor in time, then solving this block $p_t \times p_t$ ``interface'' system, then solving for the $n_t = N/p_t$ time steps on each processor in parallel.  The computational cost of the decoupling is that of $2$ block solves of systems with $n_t$ blocks on each processor (one to propagate the right-hand side data and one to propagate the coupling term between the first time step on processor $i$ in a space-time patch and the first time step on processor $i+1$).
As above, we model the cost of solving a block lower bidiagonal system as in~\eqref{eq:waveform_patch_solve} by multiplying the cost of solving one of the diagonal blocks times the dimension of the system.  Thus, in cyclic reduction, each processor effectively solves three systems of size $n_t$ (two in the ``decoupling'' step, and one in the final step) and one of size $p_t$, leading to a computational cost of $3n_t + p_t \rerevise{= 3N/p_t + p_t}$ times the cost of a spatial V-cycle over one element in time for the finest-grid system, decreasing by a factor of 1/4 for each successive coarse-grid system.
For these solves, the computational time model from~\eqref{eq:solve_cost_model} retains the first term but again replaces patch solves over $N$ temporal elements with the $3n_t + p_t$ patch solves from the cyclic reduction algorithm.  We note, however, that in the evaluation of the cost model below, we make the additional assumption that the number of V-cycles for the spatial multigrid and WRMG solvers are different, as typically fewer multigrid cycles are required to solve the system for each timestep within timestepping than are required for a single all-at-once solution of the space-time system, as with WRMG.
}

\rerevise{
  Adding these costs gives us a final cost model for timestepping for the Navier-Stokes equations as
  \begin{align}
    T_t = & T_{t,\text{comm}} + T_{t,\text{comp}} \label{eq:model_timestepping} \\
    = & 16NV_t\logtwo{M_x}\alpha + (32NV_tM_x/p_x)(q+1)(9k^2+15k+15)\beta \notag \\
    & + c_1 \left((q+1)\left(\text{spatial DoFs per core}\right)\right)^{p_1} \notag \\ & + c_2 \frac{16 V_t N M_x^2}{3p_x^2}(q+1)^{p_2}(12k^4 + 96k^3 + 240k^2 + 252k + 128)^{p_3}, \notag
  \end{align}
  where the number of spatial DoFs per core is
  \[
\frac{3(M_x-1)^2k^2 +\left(4(M_x-1)^2+6(M_x-1)\right)k +2(M_x-1)^2+4(M_x-1)+3}{p_x^2},
\]
and we use $V_t$ to denote the average number of V-cycles per timestep.
For WRMG, we have a similar accounting, leading to
\begin{align}
    T_w = & T_{w,\text{comm}} + T_{w,\text{comp}} \label{eq:model_waveform} \\
    = & 20 NV_w\logtwo{M_x}\alpha \notag \\
    & + 16V_w\left((2M_x/p_x)(N/p_t) + (8/3)(M_x/p_x)^2\right)(q+1)(9k^2+15k+15)\beta \notag \\
    & + c_1 \left((q+1)\left(\text{spatial DoFs per core}\right)\right)^{p_1} \notag \\
    & + c_2 (3N/p_t + p_t) \frac{16 V_w M_x^2}{3p_x^2}(q+1)^{p_2}(12k^4 + 96k^3 + 240k^2 + 252k + 128)^{p_3},\notag
\end{align}
where $V_w$ is the number of WRMG V-cycles needed for convergence.
Similar expressions are obtained with consistent substitutions for the heat equation, as discussed above and in the Appendix.  We note the critical difference in both the communication and computation costs, where the $\beta$ and $c_2$ terms scale directly with $N$ for timestepping, but with $3N/p_t + p_t$ for WRMG.  When these terms are dominant in the total time, this provides the opportunity for parallel speedup when $p_t = \mathcal{O}(\sqrt{N})$.  Assuming the computational $c_2$ term is dominant, this leads to a potential speedup of
\[
(3/p_t + p_t/N)\frac{V_w}{V_t}
\]
when $p_tp_x^2$ processors are used for WRMG and $p_x^2$ processors are used for timestepping.  A similar speedup is seen when the bandwidth terms (multiplying $\beta$) are compared.
}

\section{Numerical experiments} \label{sec:numerics}

Our implementation can be found in \cite{self:code}, which utilizes the
Firedrake and PETSc libraries \cite{FiredrakeUserManual,
  petsc-user-ref}. The specific versions of the Firedrake packages used
are recorded here \cite{self:code:firedrake}.  We use UMFPACK~\cite{10.1145/992200.992206} as the linear solver for the patch systems (which are always solved in serial, with the patch ``owned'' by a single processor) and MUMPS~\cite{MUMPS:1} as the linear solver for the coarsest-grid systems.
\rerevise{Timestepping results are computed using the same framework, but with temporal discretization handled via the Irksome library~\cite{farrell2021irksome, RCKirby_SPMacLachlan_2024a, BDAndrews_2026a}.}
Throughout, we initialize our coarsest spatial mesh as a regular triangulation of the unit square domain, $\Omega = [0,1]^2$,
with $10$ quadrilateral cells in both the $x$ and $y$ directions cut into 2 triangles each from bottom-right to top-left. This base spatial
mesh is then uniformly refined $\Mref$ times, forming the hierarchy of
spatial meshes. Temporally, we extrude each mesh in our hierarchy
identically and uniformly, with $N$ elements.  Simulations are performed across
\revise{$8$ Xeon $1.7$GHz CPU cores (unless otherwise stated) with a memory limit of $256$GB} per
simulation. When a simulation cannot be run within the allocated memory, we
represent this in our experiments by writing ``OoM''.

\subsection{Heat equation} \label{sec:numerics:heat}

Before applying the WRMG method detailed in \S\ref{sec:wrmg} to the
Navier-Stokes equations, we first apply it to the heat equation, as described in
Definition \ref{def:heat} with initial condition $u_0(x) = \sin{\pi x} + \cos{2\pi y}$ and homogeneous Neumann BCs. \revise{While we are primarily interested in developing the method for more complex models, the lower memory overhead of a simple diffusion equation allows us to better explore a wider range of discretization orders, to gain insight into potential performance of a parallel implementation.} Due to the rapid decay of solutions over
time, we restrict our study to the time window $t\in[0,0.02]$ and use
$N=20$ timesteps. \revised{The solution decays by less than a
  factor of 2 over this interval.}  Since this problem is linear, we solve it with WRMG-preconditioned FGMRES with a convergence tolerance requiring a reduction in the residual norm either by a relative factor of $10^{-6}$ or to below an absolute tolerance of $10^{-6}$.  \revise{Providing insightful comparison results for a serial-in-time implementation of a parallel-in-time algorithm is difficult, since we naturally expect the approach to be slower in wall-clock time to standard timestepping or a true parallel-in-time solver.  As a weak benchmark, we provide comparison results for solving the coupled (``all-at-once'') space-time linear system once with an LU factorization run on the same parallel configuration.  While this is also an imperfect comparison (since we expect LU factorization to be slower), it emphasizes the good algorithmic scaling of the WRMG solver, and the true need for algorithms other than LU factorization to solve such systems.}  \rerevise{We also provide iteration counts and timings for classical timestepping, using the monolithic multigrid approach discussed above and examined for implicit Runge-Kutta methods in~\cite{RAbuLabdeh_etal_2022a, doi:10.1137/23M1569344}, with the same stopping criteria applied to each timestep.  We note that, in this linear case, the DG-in-time discretization used here is equivalent to a Radau IIA implicit Runge-Kutta method, so it is not surprising that the timestepping solver is effective.}  In all experiments reported below, all FGMRES-WRMG solves take only 3 iterations to reach this convergence tolerance except for the case of temporal degree $3$ and spatial
degree $2$ with $\Mref=3$, which takes $4$ iterations.  \rerevise{The timestepping solves are also very efficient, averaging between 1.25 and 2.05 multigrid-preconditioned FGMRES iterations per timestep to reach this convergence tolerance in all cases, with lower iteration counts for spatial degree $1$, and all higher spatial order cases averaging 2.05 iterations per timestep.}

Table~\ref{tab:heat:dof_nnz3} presents some information about the discretization for various temporal and spatial degrees at $\Mref=3$, where the finest grid has $2\times 80^2 = 12800$ spatial elements.  At lowest order (0 in time and 1 in space), we see about $1.3\times 10^5$ total DoFs, which makes sense for a problem with $81^2 = 6561$ spatial degrees of freedom and 20 time steps, giving 131 thousand total DoFs.  For this problem, we see that \revise{Firedrake allocates} about 20 nonzeros per row, corresponding to the 7-point spatial stencil for the stiffness and mass matrices \revised{that is symbolically} coupled across 3 timesteps (\revised{since Firedrake allocates nonzeros based on possible (topological) connections in the mesh, not accounting for possible coincidental cancellation or otherwise zero entries}).  As we increase the order in time, we see the number of DoFs grows linearly, while the number of nonzero entries in the matrix grows somewhat more quickly. \revise{We note that while the nonzero counts are somewhat larger than necessary, due to this symbolic allocation, it is consistently so by the same factor for all orders.}  As we increase the spatial order, we see quadratic growth in the number of DoFs and faster growth in the number of nonzeros in the matrix. \revise{We note that for fixed number of spatial ($12800$) and temporal ($20$) elements, the increases in these DoF and nonzero counts with discretization order reflect the increasing size and density of the patch systems, as discussed above.}

\begin{table}
     \begin{center}
       \caption{Total number of space-time degrees of freedom (left) and number of nonzero matrix entries (right) for the discretization of the heat equation at various spatial and temporal degrees with $\Mref=3$.}\label{tab:heat:dof_nnz3}
       \begin{tabular}{|c||c|c|c||c|c|c|}
    \hline
     & \multicolumn{3}{c||}{Degrees of Freedom} & \multicolumn{3}{c|}{Nonzero matrix entries} \\
      \hline
      3 & 5.2e5 & 2.1e6 & 4.6e6  & 4.2e7 & 2.7e8 & 9.1e8 \\
      2 & 3.9e5 &1.6e6 & 3.5e6   & 2.4e7 & 1.5e8 & 5.1e8 \\
      1 & 2.6e5 & 1.0e6 & 2.3e6  & 1.1e7 & 6.9e7 & 2.3e8 \\
      0 & 1.3e5 & 5.2e5 & 1.2e6  & 2.6e6 & 1.7e7 & 5.7e7 \\
      \hhline{|=#===#===|}
      \diagbox[dir=SW]{time}{space} & 1 & 2 & 3 & 1 & 2 & 3 \\
      \hline
  \end{tabular}
  \end{center}
\end{table}

Table~\ref{tab:heat:time3} shows the corresponding times-to-solution, where we see somewhat sub-linear scaling in time-to-solution versus number of nonzero matrix entries using WRMG, with (for example) total time \revise{only increasing by $1.6\times$ with spatial degree 1 and temporal degree going from 1 to 2, even though the number of DoFs in the system increases by $1.5\times$, and the number of nonzeros increases by a factor of over 2.  As expected, FGMRES-WRMG outperforms direct solution of the monolithic space-time system in almost every instance (except at lowest order, where they take essentially the same time).}  \rerevise{Also as expected, without time parallelism, timestepping is significantly faster than both FGMRES-WRMG and direct solvers.  This is partly aided by the fact that there is no time-dependence in this linear problem, so the multigrid preconditioner can be built once and reused for all timesteps (unlike the Navier-Stokes case that follows).  Rebuilding the timestepper at each timestep increases these times dramatically, to about 95 seconds for the case of spatial and temporal degrees both equal to 3.}

\begin{table}
  \begin{center}
  \caption{Total computational time (in seconds) to solve the discretization of the heat equation at various spatial and temporal degrees with $\Mref=3$.  At left, results for direct solution (LU factorization).  At \rerevise{center}, times for FGMRES-WRMG.  \rerevise{At right, times for timestepping.}}\label{tab:heat:time3}
\revise{
  \begin{tabular}{|c||c|c|c||c|c|c||c|c|c|}
    \hline
     & \multicolumn{3}{c||}{direct solver} & \multicolumn{3}{c||}{FGMRES-WRMG}  & \multicolumn{3}{c|}{\rerevise{timestepping}}\\
      \hline
      3 & 122.9 & 1223.1 & 4910.5 & 35.3 & 187.0 & 784.6 & 5.7 & 8.5 & 15.5 \\
      2 &  59.4 &  553.9 & 1761.0 & 23.1 & 113.4 & 424.4 & 3.8 & 6.1 & 10.5 \\
      1 &  27.4 &  201.1 &  588.3 & 14.0 &  48.4 & 183.1 & 2.9 & 4.3 & 6.7 \\
      0 &   6.8 &   43.1 &  128.4 &  7.0 &  17.9 &  48.9 & 2.2 & 3.0 & 3.8 \\
      \hhline{|=#===#===#===|}
      \diagbox[dir=SW,width=40pt]{time}{space} & 1 & 2 & 3 & 1 & 2 & 3 & 1 & 2 & 3 \\
      \hline
  \end{tabular}}
  \end{center}
\end{table}

Similar data for the heat equation at $\Mref = 4$ is presented in Tables~\ref{tab:heat:dof_nnz4} and~\ref{tab:heat:time4}.  As expected, the problem sizes and number of nonzero matrix entries reported in Table~\ref{tab:heat:dof_nnz4} are consistently approximately 4 times those in Table~\ref{tab:heat:dof_nnz3} (with slight variations due to rounding and the effects of boundary conditions).  Notably, we see that WRMG is able to solve problems with up to \revise{9} million discrete DoFs and over \revise{1} billion nonzero matrix entries within the allotted memory footprint, while solving the monolithic system with LU fails for \revise{slightly smaller problems}.  Computational times in Table~\ref{tab:heat:time4} show much more extreme scaling than in Table~\ref{tab:heat:time3}.  At spatial degrees 1 and 2, FGMRES-WRMG now shows roughly linear scaling in the number of number of nonzero matrix entries for the problem, while LU solution shows roughly quadratic scaling \revise{(or worse), with over $26\times$} the computational time for time degree 3 and space degree 1 than for time degree 0 and space degree 1, even though the number of DoFs only increases by a factor of 4.  This cost, however, more-or-less mirrors the growth in the number of nonzeros in the matrix. \revise{When both solvers succeed, we consistently see faster results from FGMRES-WRMG, as expected, by factors ranging from about 2 for lowest-order in time up to about 7 for higher-order in time cases.}  \rerevise{In contrast, the total solution time required by timestepping grows more slowly here, still in a pre-asymptotic regime, with time-to-solution scaling like the number of DoFs (or better).  This highlights that any parallel-in-time method can only hope to outperform classical timestepping when problem sizes are large enough that strong scaling of spatial parallelism breaks down on the available resources.}

\begin{table}
  \begin{center}
  \caption{Total number of space-time degrees of freedom (left) and number of nonzero matrix entries (right) for the discretization of the heat equation at various spatial and temporal degrees with $\Mref=4$.  Results shown as ``OoM'' denote that the solver ran out of memory.}\label{tab:heat:dof_nnz4}
  \begin{tabular}{|c||c|c|c||c|c|c|}
    \hline
     & \multicolumn{3}{c||}{Degrees of Freedom} & \multicolumn{3}{c|}{Nonzero matrix entries} \\
      \hline
      3 & 2.1e6 & 8.2e6 & OoM  & 1.7e8 & 1.1e9 & OoM \\
      2 & 1.6e6 & 6.2e6 & OoM  & 9.4e7 & 6.2e8 & OoM \\
      1 & 1.0e6 & 4.1e6 & 9.3e6  & 4.2e7 & 2.7e8 & 9.1e8 \\
      0 & 5.2e5 & 2.1e6 & 4.6e6 & 1.0e7 & 6.8e7 & 2.3e8 \\
      \hhline{|=#===#===|}
      \diagbox[dir=SW]{time}{space} & 1 & 2 & 3 & 1 & 2 & 3 \\
      \hline
  \end{tabular}
  \end{center}
\end{table}

\begin{table}
  \begin{center}
    \caption{Total computational time (in seconds) to solve the discretization of the heat equation at various spatial and temporal degrees with $\Mref=4$.  At left, results for direct solution (LU factorization).  At \rerevise{center}, times for FGMRES-WRMG.  \rerevise{At right, times for timestepping.}  Results shown as ``OoM'' denote that the solver ran out of memory.}\label{tab:heat:time4}
    \revise{
  \begin{tabular}{|c||c|c|c||c|c|c||c|c|c|}
    \hline
     & \multicolumn{3}{c||}{direct solver} & \multicolumn{3}{c||}{FGMRES-WRMG} & \multicolumn{3}{c|}{\rerevise{timestepping}} \\
      \hline
      3 & 994.8 &    OoM &    OoM & 132.2 & 1064.0 & OoM    &  7.3 & 18.6 & 47.2\\
      2 & 504.0 & 4195.6 &    OoM &  86.3 &  633.6 & OoM    &  5.6 & 13.1 & 29.7\\
      1 & 182.3 & 1435.2 & 4730.7 &  50.8 &  375.5 & 1090.0 &  4.2 & 8.7 & 17.1\\
      0 &  37.3 &  251.9 &  786.4 &  21.4 &  148.7 & 373.4  &  3.0 & 5.4 & 8.9\\
      \hhline{|=#===#===#===|}
      \diagbox[dir=SW,width=40pt]{time}{space} & 1 & 2 & 3 & 1 & 2 & 3 & 1 & 2 & 3 \\
      \hline
  \end{tabular}}
  \end{center}
\end{table}

\subsection{Navier-Stokes: Chorin test
  problem} \label{sec:numerics:chorin} 

We now move on to consider the Chorin test problem for the Navier-Stokes equations, as described in
Definition \ref{def:chorin}. Recall that this test problem is
constructed with a known (manufactured) solution (see
\eqref{eqn:chorin}), so we can measure the error in our computed
solutions. Here, we choose to report the space-time $L^2$ error in
the velocity as a measure of the effectiveness of the solvers and
discretizations\revise{, although we note several other measures could be considered}.  Since the Navier-Stokes equations are nonlinear, we
use Newton's method to solve the nonlinear discretized problem.
\rerevise{We either apply this to the space-time solution, using WRMG-preconditioned FGMRES for the resulting linear solves, or to the timestepping problem, using multigrid-preconditioned FGMRES for the system at each timestep.}
In both cases, we use stopping tolerances based on the
norm of the nonlinear residual, requiring a relative reduction in the
norm by a factor of $10^{-8}$.  When using the Newton-Krylov-multigrid
approach, we use the Eisenstat-Walker stopping criteria for FGMRES to
optimize performance, with looser linear solve convergence tolerance
when we are far from solution~\cite{eisenstat1996choosing}.  We use
the same mesh hierarchy as described above for the heat equation, and
consider the time interval $t\in [0,0.1]$, which we divide uniformly
into 10, 20, and 40 temporal elements ($N$) in the experiments below.
All results presented in this section are for the case of Reynolds
number $R=10$, but we recall Remark~\ref{rem:reynolds} and note that
very similar results were seen for $R =
\mathcal{O}(10^2)$. \revised{We note that, with $R=10$, the simulated
  solution decays from $\mathcal{O}(1)$ to $\mathcal{O}(10^{-1})$ over
  the temporal domain.}

Table~\ref{tab:chorin:mg_errors} shows the space-time velocity error norms generated using the Newton-Krylov-multigrid solver as we vary the spatial and temporal finite-element orders (noting that we report the spatial order as the order of the pressure space, with the velocity space being one order higher), for various values of $N$, the number of timesteps.  We note that these results are consistent with the hypothesis that the accuracy of this discretization is strongly controlled by time-stepping error.  At temporal degree 0 and spatial degree 1, we see a clear scaling like $1/N$ in the reported errors, dropping to $1/N^2$ at temporal degree 1.  \revise{The reported errors for spatial degree 2 are identical to those for spatial degree 1, indicating that our discretization has dominant temporal error.} \rerevise{Detailed experiments, not reported here, support this conclusion, with the space-time error norm being below $10^{-5}$ with spatial degree 1 and temporal degree 3 on the grid with $\Mref = 1$, and below $10^{-6}$ with the same orders and $\Mref = 2$.}  Unfortunately, within the memory limitations of our computing platform, we were unable to run \rerevise{WRMG for} this problem at sufficiently high orders in time to balance the spatial discretization error for the smooth solution considered here.

\begin{table}
      \begin{center}
    \caption{$\Mref=3$ and $R=10$: Space-time $L^2$ errors in velocity for the Chorin test problem for various temporal and spatial orders and numbers of timesteps, $N$, using the FGMRES-WRMG solver for the Newton linearizations.  Results shown as ``OoM'' denote that the solver ran out of memory.} \label{tab:chorin:mg_errors}\revise{
    \begin{tabular}{|c||c|c|c||c|c|c|}
      \hline
      & \multicolumn{3}{c||}{$N$}       & \multicolumn{3}{c|}{$N$} \\
      \hline
      & 10 & 20 & 40 & 10 & 20 & 40 \\
      \hhline{|=#===#===|}
      1 &1.77e-3 & 4.60e-5 & 1.18e-5 & 1.77e-3 & OoM & OoM \\
      0 &6.51e-3 & 3.32e-3 & 1.70e-3 & 6.51e-3 & 3.32e-3 & 1.70e-3 \\
      \hhline{|=#===#===|}
      \diagbox[dir=SW]{time}{space} & \multicolumn{3}{c||}{1}       & \multicolumn{3}{c|}{2} \\
      \hline
    \end{tabular}}
    \end{center}
\end{table}

Time-to-solution for the \rerevise{WRMG} Newton-Krylov-multigrid solver is shown in Table~\ref{tab:chorin:wrmg_time}, \revise{with iteration counts shown in Table~\ref{tab:chorin:wrmg_iterations}.  Newton iteration counts are generally stable for this problem, with most runs requiring 4 iterations, and 5 iterations needed for both cases with spatial degree 2 and $N = 10$.  Total linear iteration counts are reasonably consistent for spatial degree 1, but both larger and more variable for spatial degree 2.  We observe that much of this variation is incurred in the final Newton step, where the conditioning of the piecewise quadratic discretization matrix appears to make it difficult to achieve the required residual reduction\rerevise{, but larger tolerances led to cases where larger space-time errors were observed}.  The wall-clock time scaling is somewhat worse than expected from the number of DoFs or nonzeros in the matrix, shown in Tables~\ref{tab:chorin:dofs} and~\ref{tab:chorin:nnz}, respectively, with increases by a factor of roughly 4 going from $N=10$ to $N=20$, and 6 or 7 going from $N=20$ to $N=40$.  We suspect this is related to memory bandwidth on the machine used for these tests, but have not confirmed the root cause.}

\begin{table}
  \begin{center}
    \caption{$\Mref=3$ and $R=10$: Time-to-solution (in seconds) for the Chorin test problem for various temporal and spatial orders and numbers of timesteps, $N$, using the FGMRES-WRMG solver for the Newton linearizations.  Results shown as ``OoM'' denote that the solver ran out of memory.} \label{tab:chorin:wrmg_time}\revise{
    \begin{tabular}{|c||c|c|c||c|c|c|}
      \hline
      & \multicolumn{3}{c||}{$N$}       & \multicolumn{3}{c|}{$N$} \\
      \hline
      & 10 & 20 & 40 & 10 & 20 & 40 \\
      \hhline{|=#===#===|}
      1 & 432.2 & 1971.7 & 13964.7 & 2868.7 & OoM & OoM \\
      0 & 110.0 & 395.3 & 2596.8 & 629.8 & 2326.1 & 16979.7 \\
      \hhline{|=#===#===|}
      \diagbox[dir=SW]{time}{space} & \multicolumn{3}{c||}{1}       & \multicolumn{3}{c|}{2} \\
      \hline
    \end{tabular}}
    \end{center}
\end{table}

\begin{table}
  \begin{center}\revise{
    \caption{$\Mref=3$ and $R=10$: Number of nonlinear and linear iterations needed for convergence for Newton's method using FGMRES-WRMG as the linear solver for the Chorin test problem for various temporal and spatial orders and numbers of timesteps, $N$, using the FGMRES-WRMG solver for the Newton linearizations.  Results shown as ``OoM'' denote that the solver ran out of memory.} \label{tab:chorin:wrmg_iterations}
    \begin{tabular}{|c||c|c|c||c|c|c|}
      \hline
      & \multicolumn{3}{c||}{$N$}       & \multicolumn{3}{c|}{$N$} \\
      \hline
      & 10 & 20 & 40 & 10 & 20 & 40 \\
      \hhline{|=#===#===|}
      1 & 4 (12) & 4 (10) & 4 (17) & 5 (26) & OoM & OoM \\
      0 & 4 (11) & 4 (11) & 4 (11) & 5 (29) & 4 (19) & 4 (19) \\
      \hhline{|=#===#===|}
      \diagbox[dir=SW]{time}{space} & \multicolumn{3}{c||}{1}       & \multicolumn{3}{c|}{2} \\
      \hline
    \end{tabular}}
    \end{center}
\end{table}
      
\begin{table}
  \begin{center}
    \caption{$\Mref=3$ and $R=10$: Number of space-time DoFs for the Chorin test problem for various temporal and spatial orders and numbers of timesteps, $N$.  Results shown as ``OoM'' denote that the solver ran out of memory.} \label{tab:chorin:dofs}
    \begin{tabular}{|c||c|c|c||c|c|c|}
      \hline
      & \multicolumn{3}{c||}{$N$}       & \multicolumn{3}{c|}{$N$} \\
      \hline
      & 10 & 20 & 40 & 10 & 20 & 40 \\
      \hhline{|=#===#===|}
      1 & 1.2e6 & 2.3e6 & 4.7e6 & 2.8e6 & OoM & OoM\\
      0 & 5.8e5 & 1.2e6 & 2.3e6 & 1.4e6 & 2.8e6 & 5.7e6 \\
      \hhline{|=#===#===|}
      \diagbox[dir=SW]{time}{space} & \multicolumn{3}{c||}{1}       & \multicolumn{3}{c|}{2} \\
      \hline
    \end{tabular}
    \end{center}
\end{table}

\begin{table}
  \begin{center}
    \caption{$\Mref=3$ and $R=10$: Number of nonzero matrix entries for the Chorin test problem for various temporal and spatial orders and numbers of timesteps, $N$.  Results shown as ``OoM'' denote that the solver ran out of memory.} \label{tab:chorin:nnz}\revise{
    \begin{tabular}{|c||c|c|c||c|c|c|}
      \hline
      & \multicolumn{3}{c||}{$N$}       & \multicolumn{3}{c|}{$N$} \\
      \hline
      & 10 & 20 & 40 & 10 & 20 & 40 \\
      \hhline{|=#===#===|}
      1 & 1.5e8 & 3.1e8 & 6.4e8 & 5.3e8 & OoM & OoM\\
      0 & 3.8e7 & 7.9e7 & 1.6e8 & 1.3e8 & 2.7e8 & 5.5e8 \\
      \hhline{|=#===#===|}
      \diagbox[dir=SW]{time}{space} & \multicolumn{3}{c||}{1}       & \multicolumn{3}{c|}{2} \\
      \hline
    \end{tabular}}
    \end{center}
\end{table}

\rerevise{For comparison, Table~\ref{tab:chorin:timestepping_time} shows time-to-solution for the same problems using a timestepping approach and the monolithic multigrid method described above to solve the resulting linearizations.  Comparing these times to those in Table~\ref{tab:chorin:wrmg_time}, we note that, even without temporal parallelism, WRMG shows a faster time-to-solution for spatial degree $1$ when $N=10$, where the memory requirements of WRMG are much below the capacity of the system, and WRMG-FGMRES requires the fewest total iterations to convergence.  As $N$ grows or we increase the spatial degree, however, we see the expected improvement in runtime with timestepping when we do not use temporal parallelism within WRMG.  Table~\ref{tab:chorin:timestepping_iterations} shows generally consistent nonlinear and linear iteration counts (now averaged over each timestep), with some improvement as $N$ increases, since the initial guess for the timestepping solution improves when the timesteps get shorter.  Within each discretization order, we see time scaling roughly consistent with the number of timesteps times the total number of linear iterations (as in the cost model proposed in \cref{sec:wrmg:performancemodel}), while costs increase roughly as expected when we change either the temporal or spatial order.
}

\begin{table}
  \begin{center}
    \caption{\rerevise{$\Mref=3$ and $R=10$: Time-to-solution (in seconds) for the Chorin test problem for various temporal and spatial orders and numbers of timesteps, $N$, using the multigrid-preconditioned FGMRES solver for the Newton linearizations within timestepping.}} \label{tab:chorin:timestepping_time}  \rerevise{
    \begin{tabular}{|c||c|c|c||c|c|c|}
      \hline
      & \multicolumn{3}{c||}{$N$}       & \multicolumn{3}{c|}{$N$} \\
      \hline
      & 10 & 20 & 40 & 10 & 20 & 40 \\
      \hhline{|=#===#===|}
      1 & 432.7 & 838.8 & 1040.7 & 772.8 & 1358.3 & 2541.9 \\
      0 & 277.2 & 533.9 & 772.0 & 395.2 & 665.6 & 1216.1 \\
      \hhline{|=#===#===|}
      \diagbox[dir=SW]{time}{space} & \multicolumn{3}{c||}{1}       & \multicolumn{3}{c|}{2} \\
      \hline
    \end{tabular}}
    \end{center}
\end{table}

\begin{table}
  \begin{center}
    \caption{\rerevise{$\Mref=3$ and $R=10$: Average number of nonlinear and linear iterations per timestep needed for convergence for Newton's method using multigrid-preconditioned FGMRES as the linear solver for the Chorin test problem for various temporal and spatial orders and numbers of timesteps, $N$, using timestepping.}} \label{tab:chorin:timestepping_iterations}  \rerevise{
    \begin{tabular}{|c||c|c|c||c|c|c|}
      \hline
      & \multicolumn{3}{c||}{$N$}       & \multicolumn{3}{c|}{$N$} \\
      \hline
      & 10 & 20 & 40 & 10 & 20 & 40 \\
      \hhline{|=#===#===|}
      1 & 3.0 (13.0) & 2.9 (11.8) & 2.3 (6.7) & 3.4 (10.8) & 3.1 (10.1) & 2.9 (10.0) \\
      0 & 3.0 (10.2) & 3.0 (9.6) & 2.6 (5.8) & 3.5 (12.3) & 3.2 (10.0) & 3.0 (9.2) \\
      \hhline{|=#===#===|}
      \diagbox[dir=SW,width=40pt]{time}{space} & \multicolumn{3}{c||}{1}       & \multicolumn{3}{c|}{2} \\
      \hline
    \end{tabular}}
    \end{center}
\end{table}

\rerevise{
  One important change between these timings and those reported for the heat equation in~\cref{tab:heat:time3,tab:heat:time4} is the substantial increase in times for timestepping on the Navier-Stokes system.  There are several factors here.  First of all, as detailed above, we have substantially more DoFs in the Navier-Stokes system than the heat equation, as well as more DoFs and nonzeros in each patch system.  Thus, while the space-time grid used in these calculations is the same as that used in the heat equation example with $\Mref=3$ and $N=20$, there is substantially more work to do in each matrix-vector product and application of the monolithic multigrid preconditioner.  Secondly, as seen in Table~\ref{tab:chorin:timestepping_iterations}, substantially more iterations are needed to reach convergence here, in comparison to the roughly 2 iterations per timestep needed for the heat equation.  Even accounting for these, however, there is a substantially larger increase in time-to-solution from the 2-4 seconds reported for the heat equation with $\Mref=3$ in Table~\ref{tab:heat:time3} and the 500-1400 seconds reported here for $N=20$.  We attribute some of this further increase in cost to the overhead of dealing with a nonlinear system (where both the Jacobian and the preconditioner must be rebuilt at each nonlinear iteration) and increased cost in solving both the patch systems and the coarsest-grid Jacobian with the larger and denser linearized systems for the mixed finite-element discretization considered here.
}

\subsection{Navier-Stokes: Lid-driven cavity} \label{sec:numerics:lid}

\begin{figure}[h!]

  \centering
  \includegraphics[width=0.3\textwidth]{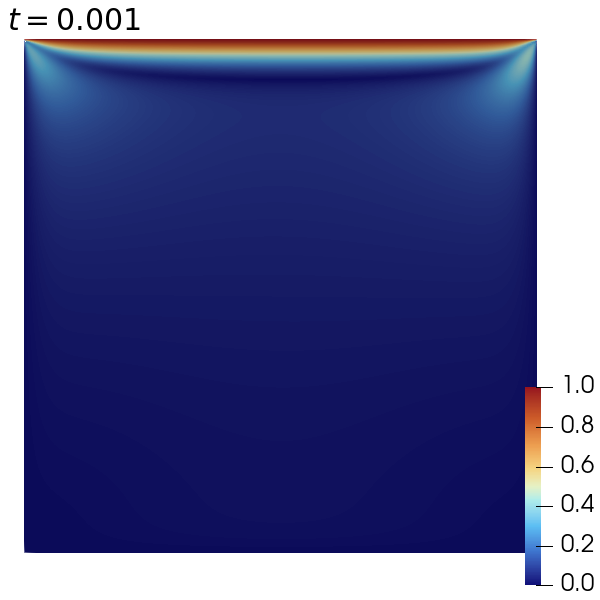}
  \includegraphics[width=0.3\textwidth]{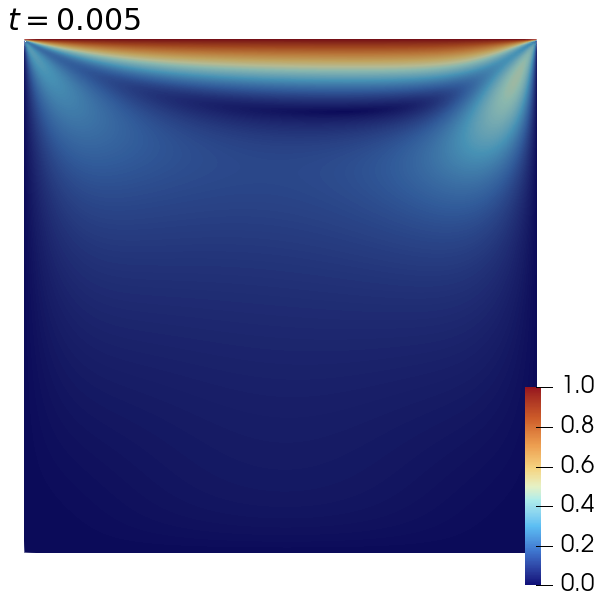}
  \includegraphics[width=0.3\textwidth]{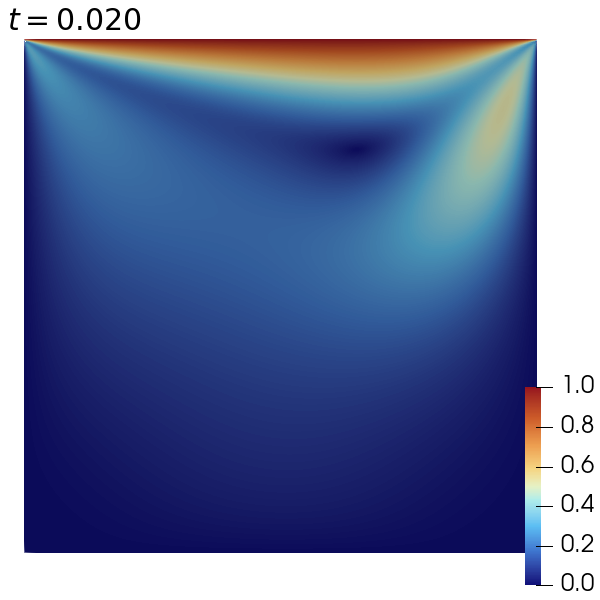}

  \caption{\rerevise{
      The magnitude of the velocity of the lid-driven cavity with $R =
      100$, $\Mref=4$, temporal degree $1$ and spatial degree $1$
      plotted at various fixed time-levels over space. We
      observe the solution changes non-trivially over our simulation
      domain.
    }   \label{fig:lidsol}}
  
\end{figure}

For our final test problem, we consider the lid-driven cavity problem for Reynolds numbers $R= 1$, $10$, and $100$.  All solver parameters are chosen as in
\S\ref{sec:numerics:chorin}, and we consider the temporal domain 
$t\in[0,0.02]$, subdivided into $N=20$ temporal elements.  We note
that the number of DoFs and nonzero matrix entries for this system at
$\Mref = 3$ are the same as those reported for $N=20$ in
Tables~\ref{tab:chorin:dofs} and~\ref{tab:chorin:nnz}, and those for
$\Mref = 2$ and $\Mref = 4$ are simply one-fourth and four times the
size, respectively, so we do not report them here. \rerevise{To
  validate the non-triviality of our solution, we visualise the
  solution at fixed time-levels in Figure \ref{fig:lidsol}. We note
  that for all Reynolds numbers considered, the dynamics remain
  non-trivial.}

We first demonstrate the effectiveness of our solvers in terms of the number of Newton steps to convergence and the total number of FGMRES-WRMG iterations required for those steps, as reported in Table~\ref{tab:lid:iterations}.  Here, we see there is some growth in the number of nonlinear and linear iterations with Reynolds number \revise{at spatial degree 1}, but the growth is quite minor.  \revise{At spatial degree 2, there is notably more growth, again likely due to conditioning issues.}  In comparison, using direct solvers to solve the \rerevise{space-time} linear systems consistently requires 2 Newton iterations for $R=1$, 3 Newton iterations for $R=10$, and 4 Newton iterations for $R=100$.  \rerevise{Table~\ref{tab:lid:timestepping_iterations} gives corresponding iteration counts, averaged per timestep, for timestepping, showing similar growth with Reynolds number, but relatively steady linear and nonlinear iteration counts as we change the order of the discretization.}
\begin{table}
    \begin{center}
    \caption{Number of nonlinear and linear iterations needed for convergence for Newton's method using FGMRES-WRMG as the linear solver for the lid-driven cavity problem for various temporal and spatial orders with varying Reynolds number, $R$, and $\Mref$.  Results shown as ``OoM'' denote that the solver ran out of memory.} \label{tab:lid:iterations}\revise{
    \begin{tabular}{|c|c||c|c|c||c|c|c|}
      \hline
      & & \multicolumn{3}{c||}{$R$}       & \multicolumn{3}{c|}{$R$} \\
      \hline
      & & 1 & 10 & 100 & 1 & 10 & 100 \\
      \hhline{|==#===#===|}
      \multirow{2}{*}{$\Mref=2$} &
      1 & 4 (10) & 4 (10) & 4 (13) & 4 (15) & 3 (19) & 4 (20) \\ 
    & 0 & 4 (10) & 4 (11) & 4 (15) & 4 (22) & 4 (17) & 4 (56) \\ 
      \hhline{|==#===#===|}
      \multirow{2}{*}{$\Mref=3$} &
      1 & 4 (8) & 4 (8) & 4 (13) & OoM & OoM & OoM \\
    & 0 & 4 (8) & 4 (8) & 4 (12) & 4 (14) & 5 (36) & 4 (33) \\ 
      \hhline{|==#===#===|}
     & \diagbox[dir=SW]{time}{space} & \multicolumn{3}{c||}{1}       & \multicolumn{3}{c|}{2} \\
      \hline
    \end{tabular}}
    \end{center}
\end{table}
\begin{table}
    \begin{center}
    \caption{\rerevise{Average number of nonlinear and linear iterations per timestep needed for convergence for Newton's method using multigrid-preconditioned FGMRES as the linear solver for the lid-driven cavity problem using timestepping, for various temporal and spatial orders with varying Reynolds number, $R$, and $\Mref$.}} \label{tab:lid:timestepping_iterations}\rerevise{
    \begin{tabular}{|c|c||c|c|c||c|c|c|}
      \hline
      & & \multicolumn{3}{c||}{$R$}       & \multicolumn{3}{c|}{$R$} \\
      \hline
      & & 1 & 10 & 100 & 1 & 10 & 100 \\
      \hhline{|==#===#===|}
      \multirow{2}{*}{$\Mref=2$} &
      1 & 1.25 (3.9) & 1.4 (4.2) & 1.9 (5.9) & 1.3 (2.7) & 1.2 (2.6) & 1.3 (3.7) \\
    & 0 & 1.35 (3.2) & 1.2 (2.7) & 1.5 (3.5) & 1.4 (2.3) & 1.2 (2.3) & 1.4 (2.8) \\
      \hhline{|==#===#===|}
      \multirow{2}{*}{$\Mref=3$} &
      1 & 1.5 (3.5) & 1.6 (3.7) & 1.6 (4.5) & 1.2 (2.4) & 1.2 (2.4) & 1.5 (4.0) \\
    & 0 & 1.4 (2.1) & 1.2 (2.4) & 1.4 (2.9) & 1.3 (2.0) & 1.1 (2.2) & 1.3 (2.7) \\
      \hhline{|==#===#===|}
     & \diagbox[dir=SW,width=40pt]{time}{space} & \multicolumn{3}{c||}{1}       & \multicolumn{3}{c|}{2} \\
      \hline
    \end{tabular}}
    \end{center}
\end{table}

Times-to-solution for both solvers, \rerevise{using FGMRES-WRMG to solve the space-time Newton linearizations or monolithic-multigrid preconditioned FGMRES to solve the linearizations at each timestep}, are reported in Tables~\ref{tab:lid:mg_times} \rerevise{and~\ref{tab:lid:timestepping_times}}, respectively. \rerevise{In both cases, times reported are \revise{mostly} independent of the Reynolds number.  In this setting, where timestepping required somewhat fewer iterations to satisfy the Newton convergence tolerance than for the Chorin test case above, we see that timestepping is between 5 and 20 times faster than the monolithic solution using WRMG-FGMRES.}

\begin{table}
  \begin{center}
    \caption{Time to solution (in seconds) needed for convergence for Newton's method using FGMRES-WRMG as linear solvers for the lid-driven cavity problem for various temporal and spatial orders with varying Reynolds number, $R$, and $\Mref$.  Results shown as ``OoM'' denote that the solver ran out of memory.} \label{tab:lid:mg_times}\revise{
    \begin{tabular}{|c|c||c|c|c||c|c|c|}
      \hline
      & & \multicolumn{3}{c||}{$R$}       & \multicolumn{3}{c|}{$R$} \\
      \hline
      & & 1 & 10 & 100 & 1 & 10 & 100 \\
      \hhline{|==#===#===|}
      \multirow{2}{*}{$\Mref=2$} &
      1 & 432.8 & 434.2 & 456.6 & 2212.5 & 1852.6 & 2344.1 \\
    & 0 & 90.6 & 90.3 & 98.8 & 468.0 & 428.6 & 627.0 \\
      \hhline{|==#===#===|}
      \multirow{2}{*}{$\Mref=3$} &
      1 & 1798.1 & 1782.4 & 1894.1 & OoM & OoM & OoM \\
    & 0 & 347.7 & 348.9 & 381.9 & 1843.8 & 2721.8 & 2317.7 \\
      \hhline{|==#===#===|}
     & \diagbox[dir=SW]{time}{space} & \multicolumn{3}{c||}{1}       & \multicolumn{3}{c|}{2} \\
      \hline
    \end{tabular}}
    \end{center}
\end{table}
\begin{table}
  \begin{center}
    \caption{\rerevise{Time to solution (in seconds) needed for convergence for timestepping using  Newton's method with monolithic-multigrid preconditioned FGMRES as the linear solver for the lid-driven cavity problem for various temporal and spatial orders with varying Reynolds number, $R$, and $\Mref$.}} \label{tab:lid:timestepping_times}\rerevise{
    \begin{tabular}{|c|c||c|c|c||c|c|c|}
      \hline
      & & \multicolumn{3}{c||}{$R$}       & \multicolumn{3}{c|}{$R$} \\
      \hline
      & & 1 & 10 & 100 & 1 & 10 & 100 \\
      \hhline{|==#===#===|}
      \multirow{2}{*}{$\Mref=2$} &
      1 & 84.4 & 87.3 & 123.6 & 119.2 & 112.7 & 133.5 \\
    & 0 & 57.2 & 49.2 & 59.9 & 61.8 & 51.7 & 62.3 \\
      \hhline{|==#===#===|}
      \multirow{2}{*}{$\Mref=3$} &
      1 & 320.3 & 324.0 & 358.4 & 422.5 & 438.5 & 591.5 \\
    & 0 & 175.9 & 167.8 & 206.4 & 196.0 & 184.5 & 214.8 \\
      \hhline{|==#===#===|}
     & \diagbox[dir=SW]{time}{space} & \multicolumn{3}{c||}{1}       & \multicolumn{3}{c|}{2} \\
      \hline
    \end{tabular}}
    \end{center}
\end{table}

\begin{figure}

  \includegraphics[width=0.9\textwidth]{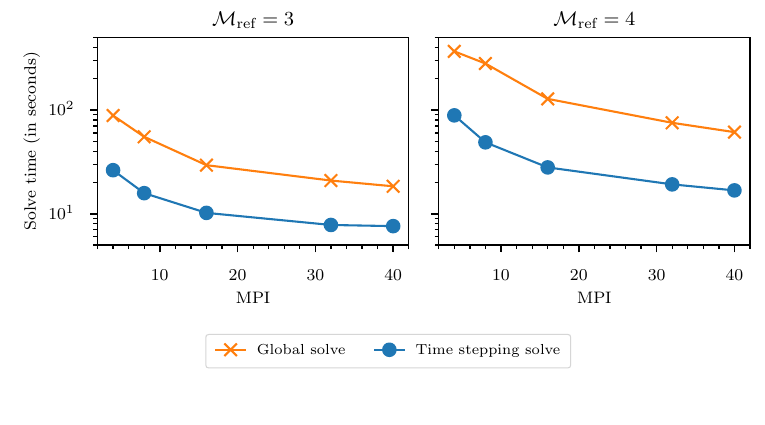}

  \caption{Time-to-solution (in seconds) for temporal degree $0$ and
    spatial degree $1$ for a simulation of the lid-driven cavity
    Navier-Stokes test case with $R=10$ for differing numbers of MPI
    ranks. Here, we compare solving the global space-time system using WRMG as a preconditioner for GMRES
    against localizing the system to a single timestep and solving
    sequentially over the timesteps.}
  \label{fig:timings}
  
\end{figure}

Finally, we consider a strong scaling experiment with increasing number of MPI ranks, comparing our monolithic Newton-Krylov-multigrid solver against the timestepping approach.  Here, we use a different machine, with 48 Xeon 2.4GHz CPU cores, to test the spatial parallelism.  \rerevise{As above,} we note that the timestepping approach results in a significantly smaller system to solve (with 1/20th the number of DoFs, in the example considered here) for each timestep, but one that must be solved 20 times to timestep to the same final solution time.  Figure~\ref{fig:timings} shows timing results for the problem with temporal degree $0$ and spatial degree $1$ at $\Mref=3$ and $\Mref=4$.  Here, we see that the timestepping solver is consistently about 5 times faster than the monolithic solve using WRMG, \rerevise{consistent with the results reported in the tables above}.  \revise{As we discuss below, this is \rerevise{also} consistent with the performance model, in that we have yet to reach the regime where the strong scaling of the spatial multigrid solver breaks down, so we are not yet where we expect to see performance benefits from WRMG.}

Clearly, in this regime (few timesteps in comparison to the spatial mesh sizes \revised{and low levels of parallelism}), the monolithic Newton-Krylov-multigrid approach is not competitive with the timestepping approach.  We note, however, that there are advantages for the monolithic space-time approach that we have not yet exploited in the results presented here, due to limitations on our testing platform.  Primarily, the WRMG approach is able to exploit finer-grained parallelism by using cyclic-reduction (or similar parallel direct solvers) to make effective use of more parallel cores by decomposing the space-time solves within the waveform relaxation~\cite{RDFalgout_etal_2015a}\revised{, as discussed above}.  Thus, on a system where many more cores are available, we can use a ``two-level'' approach to parallelism for WRMG, first decomposing the spatial grid to a sensible parallel distribution, then decomposing the temporal grid in a space-time parallel manner (that is not currently supported by Firedrake).  This also would allow us to make use of a common feature of modern parallel architectures, where the amount of memory available scales with the number of cores, supporting the higher memory footprint of monolithic space-time solvers.  We leave implementing and investigating this option for future work\revise{, but take first steps towards understanding this in the performance model considered next}.

\revised{
  \subsection{Fitting the performance model}
  \label{sec:numerics:performancemodel}

  Because we do not yet have a parallel-in-time implementation of the WRMG solver, implementing cyclic reduction or a similar algorithm to solve the space-time patch systems, we now evaluate the parallel performance model developed in~\cref{sec:wrmg:performancemodel} to get insight into the potential performance of an efficient implementation on modern HPC hardware.

  \rerevise{In this section, we use the cost models in~\cref{eq:model_timestepping,eq:model_waveform} to predict performance for a simulation with $1000$ time steps and a $160\times 160$ spatial grid with both low-order (temporal order $q=0$, spatial order $k=2$) and high-order (temporal order $q=2$, spatial order $k=4$) discretizations of the Navier-Stokes equations.  We use the parameters found above, from fitting the heat equation experiments, along with observed data on iteration counts for both the timestepping and parallel-in-time algorithms.  The parameters used in this section are summarized in \cref{tab:model_parameters}.  In particular, we assume that the parallel-in-time algorithm requires a total of $12$ multigrid
    iterations (over several Newton iterations), whereas we assume only a single Newton iteration is needed for each timestep in the timestepping approach, requiring $3$ multigrid iterations to reach convergence, roughly consistent with the data presented in~\cref{tab:lid:timestepping_iterations}.  We note that the lower timestepping iteration counts are because the solution from the previous timestep is generally expected to be a good initial guess for the current timestep, but that this is strongly impacted by the size of the timestep, and taking larger timesteps would lead to larger iteration counts.
}
    \begin{table}
      \begin{center}
        \rerevise{
    \caption{Parameters used in the cost model.} \label{tab:model_parameters}
    \begin{tabular}{|c|c|c|}
      \hline
      Symbol & Value & Explanation \\
      \hline \hline
      $N$ & 1000 & Number of timesteps \\
      $M_x$ & 160 & Elements in one dimension of spatial mesh \\
      \hline
      $V_t$ & 3 & Average V-cycles per timestep in timestepping \\
      $V_w$ & 12 & Total V-cycles to solution in WRMG \\
      \hline
      $\alpha$ & $1.5\times 10^{-3}$ & Latency parameter (seconds) \\
      $\beta$ & $1.5\times 10^{-6}$ & Inverse bandwidth (seconds per value)\\
      $c_1$ & $5.9\times 10^{-4}$ &  Interpolation operator assembly rate (seconds per DoF)\\
      $c_2$ & $ 3.5\times 10^{-7}$ & FLOPS rate for patch solves\\
      \hline
      $p_1$ & $1.29$ & Exponent in interpolation operator assembly \\
      $p_2$ & $2.45$ & Temporal DoFs exponent in patch solves \\
      $p_3$ & $1.00$ & Spatial nonzeros exponent in patch solves \\
      \hline
    \end{tabular}}
    \end{center}
      \end{table}

  Figure~\ref{fig:pm} shows how predicted performance varies with allocation of cores to spatial and temporal parallelism for the FGMRES-WRMG algorithm for two scenarios, one with ``low polynomial order'', using piecewise constants in time and the $\cpoly{3}-\cpoly{2}$ discretization in space, and one with ``high polynomial order'', using piecewise quadratics in time and the $\cpoly{5}-\cpoly{4}$ discretization in space.  In both cases, we see that using low numbers of temporal cores leads to diminishing returns when using large amounts of spatial parallelism.  Similarly, once we have ``enough'' spatial parallelism that decreasing strong scaling performance becomes apparent, we see substantial improvements in time-to-solution possible using temporal parallelism.
    \begin{figure}[h!]
      \includegraphics[width=0.9\textwidth]{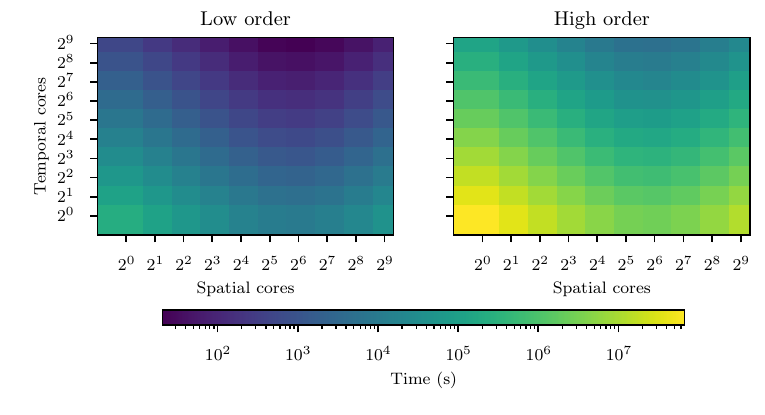}
  \caption{\revised{Time-to-solution (in seconds) as we vary temporal and
        spatial cores for the parallel-in-time performance model \rerevise{given in~\cref{eq:model_waveform} with parameters from~\cref{tab:model_parameters}} with low polynomial
        order (temporally $0$ and spatially $2$) and high
        order (temporally $2$ and spatially $4$).}}
    \label{fig:pm}
  \end{figure}

An important question raised by these results is how one should best allocate parallelism when one wants to solve a fixed problem with a fixed number of cores.  We investigate this in Figure~\ref{fig:varying}, where we consider the same low-order and high-order scenarios, but note that we consider much larger numbers of cores for the high-order case.  For both cases, we clearly see the extra work required by the WRMG algorithm that requires sufficient additional available parallelism to lead to faster runtimes.  For the low-order scenario considered here, we see decreasing strong scaling between 10 and 100 cores for simple timestepping.  When allocating only a small number of spatial cores, we see improved strong scaling up to a point, but the cost of cyclic reduction, which scales like $3N/p_t + p_t$, becomes prohibitive when $p_t$ becomes large.  Increasing the amount of spatial parallelism increases the number of total cores that can be effectively used before this increase in cost leads to poor scaling. However, with large-enough numbers of spatial cores (e.g., 1024 spatial cores), we can see real improvement in modelled wall-clock time-to-solution of about $40\times$.  A similar story is seen in the high-order case, where we see a modelled speedup of about $12\times$, when using 2M cores, with 131K spatial cores.
  \begin{figure}[h!]
    \includegraphics[width=0.9\textwidth]{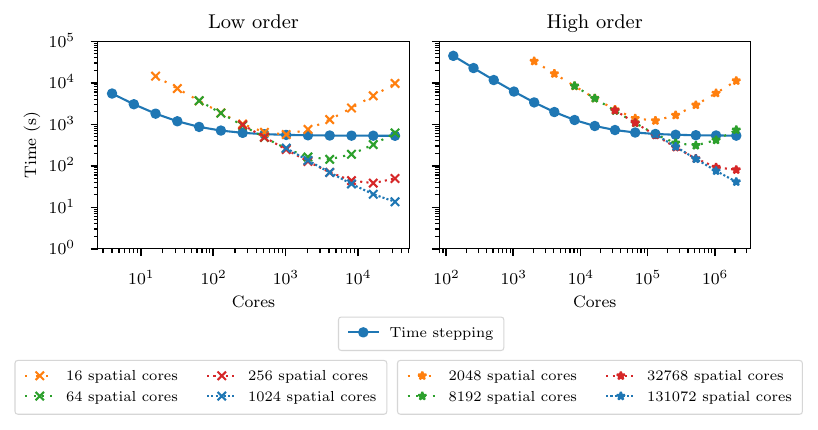}
    \caption{\revised{Time-to-solution (in seconds) for low polynomial
        order (temporally $0$ and spatially $2$) and high polynomial
        order (temporally $2$ and spatially $4$). Here we either
        consider a time-stepping approach\rerevise{, with cost modelled by~\cref{eq:model_timestepping},} and parallelise across cores
        in space, or solve globally\rerevise{, with cost modelled by~\cref{eq:model_waveform},} and parallelise in space time. \rerevise{In both cases, we use the parameters given in~\cref{tab:model_parameters}.}  We
        fix the number of spatial cores\rerevise{, $p_x^2$, while varying the total core
        count for the waveform relaxation results, and vary the total core count, $p_x^2$, for the timestepping results.}}
      \label{fig:varying}}

    \end{figure}

}

\section{Conclusion}\label{sec:conclusions}

\revise{
  In this paper, we have demonstrated that the spatial multigrid relaxation scheme proposed for the Taylor-Hood discretization of the Navier-Stokes equations in~\cite{Rafiei_2024} extends to an effective monolithic Newton-Krylov-multigrid solver for space-time finite-element discretizations.}  The solver is shown to be effective for two-dimensional flow problems, both for the simple Chorin vortex problem and the more complicated lid-driven cavity problem.  \revise{Small-scale (serial-in-time) performance results show promise for the methodology, which is further supported by the parallel performance model, which suggests substantial speedups may be possible given sufficient parallel resources.}

Due to a software \revise{and memory} limitations, we are unable to extend the method to three-dimensional flows at this time.  Overcoming these limitations is a key step in future work.  Also interesting is the extension of the method to a true space-time parallel algorithm, using variants of the cyclic reduction methodology as examined, for example, in~\cite{RDFalgout_etal_2015a}, to parallelize the solves within waveform relaxation.  This would allow a fair comparison between the method investigated herein and standard parallel-in-space timestepping in the common paradigm where there are many more cores available than can be effectively used when strong scaling a fixed spatial problem.  \revise{This would also allow comparison to other parallel-in-time approaches for the same discretizations.  Nonetheless, a parallel performance model is developed to show the potential speedups of the WRMG solver, which we believe justifies the work required to develop an efficient implementation.}  Finally, we note that the use of monolithic space-time discretizations, as considered here, allows the possibility of both local-in-space and local-in-time adaptivity, and that extensions of waveform relaxation multigrid (and other efficient monolithic solvers) to this paradigm, where classical timestepping is not directly feasible, would be of great interest.

\appendix

\rerevise{
\section{Details of parallel performance model}

In this Appendix, we provide some detailed calculations to support the discussion of the parallel performance model introduced in~\cref{sec:wrmg:performancemodel}.

\subsection{Nonzero entries in the Navier-Stokes patch matrices}
\label{ssec:patch_nonzeros}

We count the number of nonzero entries in each patch matrix by considering the number of nonzeros in each row of the patch matrix, focusing on those rows associated with velocity DoFs.  By separately the nonzero connections between velocity DoFs (that appear once in the patch matrix) and between velocity and pressure DoFs (that appear twice, once in the (1,2) block of the patch matrix and once in the (2,1) block), we arrive at an expression for the full number of nonzeros in each patch matrix.
\begin{itemize}
\item The central vertex has 2 velocity DoFs, and nonzeros entries in the matrix appear for its connections to the 2 velocity DoFs at all 7 vertices in the patch, all 12 edges in the patch, and all 6 elements in the patch.  There are also nonzeros associated with these velocity DoFs and the central vertex pressure DoF, the pressure DoFs along the 6 interior edges, and the pressure DoFs interior to the 6 elements on the patch.
\item The 6 boundary vertices all have 2 velocity DoFs, and nonzero entries in the matrix appear for their connections to the 2 velocity DoFs on each of 4 vertices, all of the edge velocity DoFs on 5 edges, and all of the interior element velocity DoFs on 2 elements.  There are also nonzeros for these velocity DoFs with the central vertex pressure DoF, pressure DoFs on 3 edges, and pressure DoFs interior to 2 elements.
\item All of the velocity DoFs on the 6 ``interior'' edges to the patch (those adjacent to the central vertex) have nonzero connections to the 2 velocity DoFs on each of 4 vertices, all of the edge velocity DoFs on 5 edges, and all of the interior element velocity DoFs on 2 elements.  There are also nonzeros for these velocity DoFs with the central vertex pressure DoF, pressure DoFs on 3 edges, and pressure DoFs interior to 2 elements.
\item All of the velocity DoFs on the 6 ``boundary'' edges of the patch have nonzero connections to the 2 velocity DoFs on each of 3 vertices, all of the edge velocity DoFs on 3 edges, and all of the velocity DoFs interior to one element.  They also have nonzero connections to the central vertex pressure DoF, the pressure DoFs on 2 interior edges, and those interior to one element.
  \item All of the velocity DoFs interior to the 6 elements in the patch have nonzero connections to the 2 velocity DoFs on each of 3 vertices, all of the edge velocity DoFs on 3 edges, and all of the velocity DoFs interior to one element (itself).  They also have nonzero connections to the central vertex pressure DoF, the pressure DoFs on 2 interior edges, and those interior to one element.
\end{itemize}

To complete the calculation of number of nonzeros in the patch, we note that each component of the $\bpoly{k+1}$ velocity field has
\begin{itemize}
\item one DoF at each vertex,
\item $k$ DoFs per edge, and
\item $k(k-1)/2$ DoFs interior to each element,
\end{itemize}
while the $\bpoly{k}$ pressure fields has one DoF at each vertex, $k-1$ DoFs per edge, and $(k-1)(k-2)/2$ DoFs interior to each element.  Thus, the calculations above give us a total number of nonzeros in the velocity block of
\begin{align*}
  & 4\left(7 + 12k + 6k(k-1)/2\right) + 12\left(8 + 10k + 4k(k-1)/2\right) \\
  & + 12k\left(8  +10k + 4k(k-1)/2\right) + 12k\left(6+6k+2k(k-1)/2\right) \\
  & + 12\left(k(k-1)/2\right)\left(6+6k+2k(k-1)/2\right),
\end{align*}
and the total number of nonzeros in one of the pressure blocks as
\begin{align*}
  &2\left(1+6(k-1)+6(k-1)(k-2)/2\right) + 12\left(1+3(k-1)+2(k-1)(k-2)/2\right) \\
  &+ 12k\left(1+3(k-1)+2(k-1)(k-2)/2\right) + 12k\left(1+2(k-1)+(k-1)(k-2)/2\right) \\
  &+ 12\left(k(k-1)/2\right)\left(1+2(k-1)+(k-1)(k-2)/2\right).
\end{align*}
Simplifying these expressions gives a total number of nonzeros in a patch matrix as
\begin{align*}
  6k^4& + 60k^3 + 198k^2 + 264k + 124 + 2\left(3k^4 + 18k^3 + 21k^2 - 6k + 2\right) \\
  & = 12k^4 + 96k^3 + 240k^2 + 252k + 128.
  \end{align*}

\subsection{Parallel communication model}
\label{ssec:halos}

For spatial multigrid (timestepping), each processor (simultaneously) needs to send messages to each of its 4 neighbours to exchange updated values of the approximate solution for each relaxation sweep.  If we do $V$ V(2,2)-cycles with \rerevise{$\logtwo{M_x}$} levels to solve each timestep, this accounts for \rerevise{$4NV\logtwo{M_x}$} total relaxation sweeps, which we multiply by 4 to get the total multiplier of $\alpha$ in the communication cost.  Each of these messages contains solution values within a strip of length equal to the number of vertices along one edge (on a given level of the hierarchy) and width equal to two elements, in order that we can accurately compute residuals at all points within all patches owned by the processor.  We approximate this by the length of the edge times the number of DoFs in a spatial patch over one temporal element.  Summing this over all levels gives a multiplier of $\beta$ of $32NVM_x/p_x$ times the number of number of DoFs in a spatial patch over one temporal element, where the factor of 32 comes from multiplying 4 edges to each processor's domain times 4 relaxation sweeps per cycle times a factor of 2 from the geometric series of edge lengths.  Thus, for spatial multigrid within timestepping for the Navier-Stokes problem, we have a total communication time model of
\[
T_{t,\text{comm}} = 16NV\logtwo{M_x}\alpha + (32NVM_x/p_x)(q+1)(9k^2+15k+15)\beta,
\]
where we use the added subscript $t$ to signify timestepping, in contrast to WRMG with subscript $w$ below.
For the heat equation, the only change is in the number of spatial DoFs per patch, which is $3k^2-3k+1$ instead of $9k^2+15k+15$.

For parallel-in-time WRMG, there are two important changes in the communication model.  First, the cost of a halo exchange for each relaxation sweep now needs to include a space-time halo.  Assuming $p_t$ processors used to parallelize in the temporal direction, this multiplies the size of each message from above by $n_t = N / p_t$.  
Secondly, using cyclic reduction, we require an additional communication step to form the interface system on each relaxation sweep.  Here, we model this as an all-to-all communication, where each processor sends two ``vectors'', of message size equal to the local number of patches times the size of the space-time patch over one temporal element, so that all processors can solve the interface system for the first time-step on all processors.  Summing this over all levels gives a factor of $4/3$ times the finest-grid cost. Thus, for WRMG for the Navier-Stokes equations, we have a total communication time model of
\begin{align*}
  T_{w,\text{comm}} = & 20 NV\logtwo{M_x}\alpha\\
  & + 16V\left((2M_x/p_x)(N/p_t) + (8/3)(M_x/p_x)^2\right)(q+1)(9k^2+15k+15)\beta.
\end{align*}

\section*{Acknowledgments}

We gratefully thank two anonymous reviewers for their detailed comments that improved this work.}

\bibliographystyle{siamplain}
\bibliography{rewaves}

\end{document}